\documentclass[12pt]{article}
\usepackage{graphicx}
\usepackage{amsmath}
\usepackage{amssymb}
\usepackage{theorem}
\usepackage{enumerate} 
\usepackage{color}

    \usepackage{silence}
\usepackage{hyperref}

\numberwithin{equation}{section}

\newtheorem{thm}{Theorem}[section]
\newtheorem{lemma}[thm]{Lemma}

\newtheorem{prop}[thm]{Proposition}
\newtheorem{cor}[thm]{Corollary}
{\theorembodyfont{\rmfamily}

\newtheorem{example}[thm]{Example}

\newtheorem{rmk}[thm]{Remark}
}

\newcommand{\qed}{\hfill \mbox{\raggedright \rule{.07in}{.1in}}}
 
\newenvironment{proof}{\vspace{1ex}\noindent{\bf
Proof}\hspace{0.5em}}{\hfill\qed\vspace{1ex}}
\newenvironment{pfof}[1]{\vspace{1ex}\noindent{\bf Proof of
#1}\hspace{0.5em}}{\hfill\qed\vspace{1ex}}

\newcommand{\R}{{\mathbb R}}

\newcommand{\Z}{{\mathbb Z}}

\newcommand{\cB}{{\mathcal B}}
\newcommand{\cC}{{\mathcal C}}

\newcommand{\cD}{{\mathcal D}}
\newcommand{\cI}{{\mathcal I}}
\newcommand{\cW}{{\mathcal W}}
\renewcommand{\H}{{\mathbb H}}
\newcommand{\barH}{{\overline{\H}}}

\newcommand{\essinf}{\operatorname{ess\, inf}}
\newcommand{\infI}{{\SMALL \inf_I}}
\newcommand{\supI}{{\SMALL \sup_I}}

\newcommand{\Leb}{\operatorname{Leb}}
\newcommand{\Int}{\operatorname{Int}}
\newcommand{\Var}{\operatorname{Var}}
\newcommand{\diam}{\operatorname{diam}}
\newcommand{\sgn}{\operatorname{sgn}}

\newcommand{\supp}{\operatorname{supp}}

\newcommand{\BV}{{\rm BV}}
\renewcommand{\Re}{\operatorname{Re}}
\renewcommand{\Im}{\operatorname{Im}}

\newcommand{\eps}{{\epsilon}}

\newcommand{\SMALL}{\textstyle}

\newcommand{\vertiii}[1]{{\left\vert\kern-0.25ex\left\vert\kern-0.25ex\left\vert #1
    \right\vert\kern-0.25ex\right\vert\kern-0.25ex\right\vert}}

\title{Renewal theorems and mixing for non Markov flows with infinite measure}

\author{
Ian Melbourne \thanks{Mathematics Institute, University of Warwick, Coventry, CV4 7AL, UK}
\and
Dalia Terhesiu
\thanks{Department of Mathematics, University of Exeter, Exeter EX4 4QF, UK}
}

\date{6 March 2017; updated 19 December 2018}

\begin{document}

 \maketitle

\begin{abstract}
We obtain results on mixing for a large class of (not necessarily Markov) infinite measure semiflows and flows.  Erickson proved, amongst other things, a strong renewal theorem in the corresponding i.i.d.\ setting.
Using operator renewal theory, we extend Erickson's methods to the deterministic (i.e.\ non-i.i.d.) continuous time setting and obtain results on mixing as a consequence.

Our results apply to intermittent semiflows and flows
of Pomeau-Manneville type (both Markov and nonMarkov), and to semiflows and flows over Collet-Eckmann maps with nonintegrable roof function.

\end{abstract}

\section{Introduction}

Recently, there has been increasing interest in the investigation of mixing properties for infinite measure-preserving dynamical systems~\cite{Aaronson13,BMTapp,Gouezel11,KautzschKessebohmerSamuel16,KautzschKessebohmerSamuelStratmann15,Lenci10,Lencisub,LiveraniTerhesiu16,M15,MT12,MT17,Terhesiu15,Terhesiu16, Thaler00}.
Most of these results are for discrete time noninvertible systems.

For results on semiflows preserving an infinite measure, we refer to~\cite{MT17} (the Markov case) and~\cite{BMTapp} (which does not assume a Markov structure).
The setting is that $F:Y\to Y$ is a mixing uniformly expanding map defined on a probability space $(Y,\mu)$ and
$\tau:Y\to\R^+$ is a nonintegrable roof function with regularly varying tails:
\begin{align} \label{eq:tail}
 \mu(y\in Y:\tau(y)>t)=\ell(t)t^{-\beta}\quad\text{for various ranges of $\beta\in[0,1]$.}
\end{align}
Here,  $\ell:[0,\infty)\to[0,\infty)$ is a measurable slowly varying function (so $\lim_{t\to\infty}\ell(\lambda t)/\ell(t)=1$ for all $\lambda>0$).
Consider the suspension $(Y^\tau,\mu^\tau)$ and suspension semiflow $F_t:Y^\tau\to Y^\tau$ (the standard definitions are recalled in Section~\ref{sec:semiflow}).
The aim is to prove a mixing result of the form
\[
\lim_{t\to\infty}a_t\int_{Y^\tau}v\,w\circ F_t\,d\mu^\tau=
\int_{Y^\tau}v\,d\mu^\tau \int_{Y^\tau}w\,d\mu^\tau,
\]
for a suitable normalisation $a_t\to\infty$ and suitable classes of observables
$v,w:Y^\tau\to\R$.

Under certain hypotheses,~\cite{BMTapp,MT17} obtained results on mixing and rates of mixing for such semiflows.
The hypotheses were of two types: (i) assumptions on ``renewal operators'' associated to the transfer operator of $F$ and the roof function $\tau$, and (ii) Dolgopyat-type assumptions of the type used to obtain mixing rates for finite measure (semi)flows~\cite{Dolgopyat98b}.

As pointed out to us by Dima Dolgopyat, P\'eter N\'andori and Doma Sz\'asz, mixing for indicator functions can be regarded as a local limit theorem and hence
hypotheses of type (ii) should not be necessary.  

In this paper, we show that operator renewal-theoretic assumptions (i) are indeed sufficient for obtaining the mixing results in~\cite{BMTapp,MT17}.   The abstract framework in~\cite{BMTapp} turns out again to be flexible enough to cover nonMarkov situations.  Moreover, our main results extend to flows and we are able to treat large classes of observables $v,w$.
(Conditions of type (i) alone are not sufficient for obtaining rates of mixing; the best results remain those in~\cite{BMTapp}.)

The analogous probabilistic results go back to Erickson~\cite{Erickson70}
who obtained {\em strong renewal theorems} in an i.i.d.\ continuous time framework under the assumption $\beta\in(\frac12,1]$.  
(In the discrete time setting, see~\cite{GarsiaLamperti62} for the
i.i.d.\ case and~\cite{MT12} for the deterministic case.)
Our results on mixing when $\beta\in(\frac12,1]$ for semiflows (Corollary~\ref{cor:mixing} and the extensions in Section~\ref{sec:obs}) and for flows (Theorem~\ref{thm:flow}),
are proved by adapting Erickson's methods to the deterministic setting.

For $\beta\le\frac12$, additional hypotheses are needed on the tail of $\tau$
to obtain a strong renewal theorem (and hence mixing) even for discrete time; see~\cite{CaravennaDoneysub, Doney97, GarsiaLamperti62} for i.i.d.\ results and~\cite{Gouezel11} for deterministic results (see also~\cite{Terhesiu16} for higher order theory in both the i.i.d.\ and deterministic settings).
For the continuous time case, Dolgopyat \& N\'andori~\cite{DN-PC} obtain strong renewal theorems for a class of Markov semiflows including the range $\beta\le\frac12$ (again under extra hypotheses on the tail $\mu(\tau>t)$), though our main examples seem beyond their framework.
In the absence of additional tail hypotheses,~\cite{Erickson70} showed how to obtain a partial result in the probabilistic setting
with limit replaced by $\liminf$.   In Corollary~\ref{cor:inf},
we obtain such a $\liminf$ result for semiflows with $\beta\in(0,\frac12]$.

We now describe two families of examples to which our results apply.
For definiteness, we restrict to our main mixing result Corollary~\ref{cor:mixing} which applies when $\beta\in(\frac12,1]$.  
(Corollaries~\ref{cor:mixingav} and Corollary~\ref{cor:inf} hold for all $\beta\in(0,1]$.)

\begin{example}[NonMarkovian intermittent semiflows and flows]
\label{ex:neutral}
Consider the map $f:[0,1]\to [0,1]$ given by
\begin{align*}
f(x) = x(1+c_1 x^{1/\beta}) \bmod 1 \quad\text{where
$\beta \in \SMALL(\frac12,1],\, c_1 >0$}.
\end{align*}
This is an example of an AFN map~\cite{Zweimuller98}, namely a nonuniformly expanding one-dimensional map with at most countably (in this case finitely) many branches with finite images and satisfying Adler's distortion condition
$\sup |f''|/|f'|^2 < \infty$.
Up to scaling, there is a unique absolutely continuous invariant measure
$\mu_0$.  The measure $\mu_0$ is infinite and the density has a singularity at the neutral fixed point $0$.

Let $\tau_0:[0,1]\to[1,\infty)$ be a roof function of bounded variation and H\"older continuous,
and let $f_t$ denote the suspension semiflow on $[0,1]^{\tau_0}$ with
invariant measure $\mu_0^{\tau_0}=\mu_0\times{\rm Lebesgue}$.
Note that there is now a neutral periodic orbit of period $\tau_0(0)$.

In~\cite{BMTapp}, under a Dolgopyat-type condition on $\tau_0$ and for 
sufficiently regular observables $v$ and $w$ supported away from the neutral periodic orbit, we proved a mixing result with rates and higher order asymptotics.
Here we obtain the mixing result without requiring the Dolgopyat-type condition or high regularity for the observables.
It suffices that $f_t$ has two periodic orbits (other than the neutral periodic orbit)
whose periods have irrational ratio.
Define
$m(t)=\begin{cases} \log t & \beta=1 \\ t^{1-\beta} & \beta\in(\frac12,1)\end{cases}$.
We show that
 \begin{align} \label{eq:ex}
\SMALL \lim_{t\to\infty}m(t)\int v\,w\circ f_t\,d\mu_0^{\tau_0}=
{\rm const}\int v\,d\mu_0^{\tau_0}\int w\,d\mu_0^{\tau_0},
\end{align}
where the constant depends only on $f$ and $\tau_0$.
Here, $v$ is any continuous function supported away from the neutral periodic orbit and $w$ is any integrable function.

\begin{rmk}
If $c_1$ is a positive integer, then $f$ is Markov and is a special case of the class of maps considered by~\cite{Thaler83}.
In this case, it suffices that $\tau_0$ is H\"older continuous.
Moreover, it follows from~\cite{DN-PC} that the mixing result~\eqref{eq:ex} holds for all $\beta\le1$.  When $c_1$ is not an integer, $f$ is not Markov and~\cite{DN-PC} does not apply, as far as we can tell, regardless of the value of $\beta$.
\end{rmk}

As in~\cite{LiveraniTerhesiu16,M15}, we can also consider solenoidal flows with 
a neutral periodic orbit.   Our results on mixing apply equally to such flows, see Remark~\ref{rmk:PH2}.
\end{example}

\begin{example}[Suspensions over unimodal maps]
\label{ex:unimodal}

We consider a class of examples studied in~\cite[Example~1.2]{BMTapp}.
Under a Dolgopyat-type condition on $\tau_0$ and for 
sufficiently regular observables $v$ and $w$, we proved a mixing result with rates and higher order asymptotics.
Again, the emphasis is now on mixing rather than mixing rates, with significantly relaxed hypotheses on the roof function and the observables.

Let $f:[0,1]\to [0,1]$ be a $C^2$ unimodal map with unique non-flat critical point $x_0\in(0,1)$.
We suppose further that $f$ is {\em Collet-Eckmann}~\cite{ColletEckmann83}: there are constants $C>0$, $\lambda_{{\textsc{\tiny CE}}}>1$ such that $|(f^n)'(fx_0)|\ge C\lambda_{{\textsc{\tiny CE}}}^n$ for all $n\ge1$.
It follows~\cite{Jak81} that there is a unique acip $\mu_0$ that is mixing up to a finite cycle.  We restrict to the case when $\mu_0$ is mixing.
Finally, we suppose that $x_0$ satisfies {\em slow recurrence} in the sense
that $\lim_{n\to\infty}n^{-1}\log |f^nx_0-x_0|=0$.

Consider a roof function $\tau_0:[0,1]\to\R^+$ of the form
$\tau_0(x)=g(x)|x-x_0|^{-1/\beta}$ where $\beta\in(\frac12,1)$ and $g:[0,1]\to(1,\infty)$ is differentiable, and 
form the suspension semiflow $f_t:[0,1]^{\tau_0}\to [0,1]^{\tau_0}$.
Suppose that $f_t$ has two periodic orbits whose periods have irrational ratio.  
We obtain the mixing property~\eqref{eq:ex}
for any continuous function $v$ supported in $[0,1]\times[0,1]$ and any integrable $w$.
\end{example}

The remainder of this paper is organised as follows.
In Section~\ref{sec:results}, we describe the operator renewal-theoretic hypotheses required in this paper and we state a strong renewal theorem for $\beta\in(\frac12,1]$ as well as related results for $\beta\le\frac12$.
In Section~\ref{sec:semiflow}, we show how these results lead to mixing properties for semiflows.  
Sections~\ref{sec:steps} and~\ref{sec:complete} are devoted to the proof of the strong renewal theorem, while Sections~\ref{sec:LLT} and~\ref{sec:pfthmaver} contain the proofs of the remaining results in Section~\ref{sec:results}.
Section~\ref{sec:BMT} contains prerequisites from operator renewal theory.

Corollary~\ref{cor:mixing} (mixing for semiflows) is stated for observables that are certain indicator functions.  This restriction is relaxed considerably in Section~\ref{sec:obs}.  The corresponding result for flows is stated and proved in Section~\ref{sec:flow}.

Finally, in Section~\ref{sec:ex} we return to
Examples~\ref{ex:neutral} and~\ref{ex:unimodal}.

\vspace{-2ex}
\paragraph{Notation}
We use ``big O'' and $\ll$ notation interchangeably, writing $a_n=O(b_n)$ or $a_n\ll b_n$
if there is a constant $C>0$ such that
$a_n\le Cb_n$ for all $n\ge1$.  Also, we write $a_n\sim b_n$ if $\lim_{n\to\infty} a_n/b_n=1$.

\section{Strong renewal theorem for continuous time deterministic systems}
\label{sec:results}

Let $(Y,\mu)$ be a probability space and let
$F:Y\to Y$ be an ergodic and mixing measure-preserving transformation.
Let $\tau:Y\to\R^+$ be a measurable nonintegrable function bounded away from zero.
For convenience, we suppose that $\essinf \tau>1$.
Throughout we assume the regularly varying tail condition~\eqref{eq:tail}.

Let $\tau_n=\sum_{j=0}^{n-1}\tau\circ F^j$.
 Given measurable sets $A,B\subset Y$, define the renewal measure
\begin{align} \label{eq:U}
U_{A,B}(I)=\sum_{n=0}^\infty \mu(y\in A\cap F^{-n}B:\tau_n(y)\in I),
\end{align}
for intervals $I\subset\R$. We write
$U_{A,B}(x)=U_{A,B}([0,x])$ for $x>0$.
Our aim is to generalise~\cite[Theorems 1 and 2]{Erickson70} to this set up. 
That is, we want to obtain the asymptotics of $U_{A,B}(t+h)-U_{A,B}(t)$ for any $h>0$. 

With the same notation as in~\cite{BMTapp}, let 
$\barH=\{\Re s\ge0\}$.  Given $\delta>0$ and $L>0$, let
$\H_{\delta,L}=(\barH\cap B_\delta(0))\cup\{ib:|b|\le L\}$.
Define the family of operators for $s\in\barH$, 
\[
\hat R(s):L^1(Y)\to L^1(Y),\qquad\hat R(s)v=R(e^{-s\tau}v).
\]
Here $R:L^1(Y)\to L^1(Y)$ is the transfer operator for $F$
(so $\int_Y Rv\,w\,d\mu=\int_Y v\,w\circ F\,d\mu$ for all $v\in L^1(Y)$, $w\in L^\infty(Y)$).

We assume that there exists $p_0\ge1$, and for each $p\in(p_0,\infty)$, $\gamma\in(0,\beta)$ and $L>0$ there exists a Banach space $\cB=\cB(Y)$ containing constant functions, with norm $\|\,\|_{\cB}$,
and constants 
$\delta\in(0,L)$, $\alpha_0\in(0,1)$ and $C>0$ such that
\begin{description}

\parskip = -2pt
\item[\textbf{(H)}]
\begin{itemize}
\item[(i)]
$\cB$ is compactly embedded in $L^p$.
\item[(ii)]
$\|\hat R(s)^n v\|_{\cB}\leq C(|v|_p+\alpha_0^n \|v\|_{\cB})$
for all $s\in \H_{\delta,L}$, $v\in\cB$, $n\ge1$.
\item[(iii)]
$|R(\tau^{\gamma}v)|_p\le C\|v\|_{\cB}$ for all $v\in\cB$.
\end{itemize}
\end{description}
Also, most of our results require one of the following conditions:
\begin{description}

\parskip = -2pt
\item[\textbf{(S)}]
\begin{itemize}
\item[(i)] For all $b\in[-L,L]$, $b\neq0$, the spectrum of $\hat R(ib):\cB\to\cB$ does not contain~$1$.
\item[(ii)] For all $b\in[-L,L]$, $b\neq0$, the spectral radius of $\hat R(ib):\cB\to\cB$ is less than~$1$.
\end{itemize}
 \end{description}
Hypothesis (H) is similar to~\cite[hypothesis~(H1)]{BMTapp}.
The hypotheses in (S) are significant weakenings of~\cite[hypothesis~(H4)]{BMTapp}
and the diophantine ratio assumption used in~\cite{MT17} 
(Dolgopyat-type condition). 
The remaining hypotheses
in~\cite{BMTapp}, namely (H2) and (H3) (re-inducing),
are not required.

\begin{rmk}
(a) For ease of exposition, hypothesis~(H) is stated on the half-plane $\barH$, though we only use $s$ real and $s$ imaginary in this paper.
For Theorems~\ref{thm:main} and~\ref{thm:inf}, we can take $s=ib$, $b\in[-L,L]$ in (H)(ii).
For Theorem~\ref{thm:average}, we can take $s=a$, $a\in[0,\delta)$ in (H)(ii).
\\[.75ex]
(b) For our main results Theorem~\ref{thm:main} and Corollary~\ref{cor:mixing}, it suffices that $\gamma>1-\beta$ (this is possible since $\beta>\frac12$ in those results).  For our other results which include $\beta\le\frac12$, it suffices that $\gamma>0$.

In addition, as in~\cite{BMTapp}, there exists $p_0\ge1$ depending only on $\beta$ and~$\gamma$ such that~(H) is required to hold only for one value of $p>p_0$.
\end{rmk}

\begin{rmk}
In the simplest setting, studied in~\cite{MT17}, where the map $F:Y\to Y$ is Gibbs-Markov~\cite{Aaronson,AaronsonDenker01}, 
hypothesis~(H) is satisfied 
with $\cB$ a symbolic H\"older space and $p=\infty$.
See~\cite[Remark~2.4]{BMTapp} and~\cite[Proposition~3.5]{MT17} for further details.  This includes the case of Markovian intermittent semiflows.

As explained in Section~\ref{sec:AFU}, this situation generalizes to the case when $F$ is an AFU map (i.e.\ an AFN map as defined in Example~\ref{ex:neutral} but uniformly expanding instead of nonuniformly expanding),
with $\cB$ consisting of bounded variation functions, enabling us to treat the nonMarkovian intermittent semiflows in Example~\ref{ex:neutral}.
\end{rmk}

Define
\[
d_\beta= \begin{cases} \frac{1}{\pi}\sin\beta\pi & \beta<1  \\ \hphantom{XX} 1 & \beta=1 \end{cases}, \qquad
 m(t)= \begin{cases}
\ell(t) t^{1-\beta} & \beta<1 \\ \int_1^t \ell(s)s^{-1}\,ds &  \beta=1\end{cases}.
\]
Throughout we suppose that
$A,B\subset Y$ are measurable and that $1_A\in\cB$.  

Our main result generalizes~\cite[Theorem 1]{Erickson70} to the present non i.i.d.\ set up:

\begin{thm}[Strong renewal theorem]
\label{thm:main}	 Assume $\mu(\tau>t)=\ell(t)t^{-\beta}$ where $\beta\in(\frac12,1]$.
Suppose that (H) and (S)(i) holds.  Then for any $h>0$, 
\[
\lim_{t\to\infty}m(t)(U_{A,B}(t+h)-U_{A,B}(t))=d_{\beta}\mu(A)\mu(B)h.
\]
\end{thm}

As discussed in the introduction,
additional hypotheses are needed to obtain a strong renewal theorem when $\beta\le\frac12$.
However, generalizing~\cite[Theorem 2]{Erickson70} 
to the present non i.i.d.\ set up, we still obtain a $\liminf$ result:

\begin{thm}\label{thm:inf}	 Assume $\mu(\tau>t)=\ell(t)t^{-\beta}$ where $\beta\in(0,1)$.
Suppose that (H) and (S)(ii) holds.  Then for any $h>0$, 
\[
\liminf_{t\to\infty}m(t)(U_{A,B}(t+h)-U_{A,B}(t))=d_{\beta}\mu(A)\mu(B)h.
\]
\end{thm}

\begin{rmk} In the i.i.d.\ setting, results of this type are first due to~\cite{GarsiaLamperti62} for discrete time and $\beta<1$.  The results of~\cite{Erickson70} extended~\cite{GarsiaLamperti62} to continuous time and incorporated the case $\beta=1$.
\end{rmk}

For the proof of Theorem~\ref{thm:inf}, we will need the following result
which gives the asymptotics of $U_{A,B}$ for the entire range $\beta\in [0,1]$.
This implies a   property for the semiflow $F_t$ 
 known as {\em weak rational ergodicity} \cite{Aaronson,AdamsSilva16} 
(see Corollary~\ref{cor:mixingav} below)
and thus is of interest in its own right.

\begin{thm}\label{thm:average}
Assume $\mu(\tau>t)=\ell(t)t^{-\beta}$ where $\beta\in [0,1]$.
Suppose that (H) holds.  Then
\[
\lim_{t\to\infty}t^{-1}m(t)U_{A,B}(t)=D_{\beta}\,\mu(A)\mu(B),
\]
where $D_\beta=\{\Gamma(1-\beta)\Gamma(1+\beta)\}^{-1}$ if $\beta\in(0,1)$
and $D_0=D_1=1$.
\end{thm}

For the proof of Theorem~\ref{thm:inf}, we will also require the following local limit theorem with error term which may also be of interest in its own right.
Let 
\[
\SMALL q_\beta(t)=\frac{1}{2\pi}\int_{-\infty}^\infty e^{ibt}e^{-c_\beta|b|^\beta}\,db, \qquad
c_\beta= 
i\int_0^\infty e^{-i\sigma}\sigma^{-\beta}\,d\sigma.
\]

\begin{thm}[LLT]
\label{thm:LLT}
Assume the setting of Theorem~\ref{thm:inf} with $\beta\in(0,1)$. 
Let $d_n>0$ be an increasing sequence with
$d_n\to\infty$ such that
$n\mu(\tau>d_n)=n\ell(d_n)d_n^{-\beta}\to 1$, as $n\to\infty$.
Then for any $h>0$ there exists $e_n>0$ with $\lim_{n\to\infty}e_n=0$ such that for all $t>0$, $n\ge1$,
\[
\Big|\mu(y\in A\cap F^{-n}B:\tau_n(y)\in[t,t+h])-\frac{h}{d_n}q_\beta(t/d_n)\mu(A)\mu(B)\Big|\le\frac{e_n}{d_n}.
\]
\end{thm}

\paragraph{Alternative hypotheses}

In certain examples, such as those where $F:Y\to Y$ is modelled by a Young tower with exponential tails~\cite{Young98}, hypothesis~(H)(iii) is problematic.
In such cases, it is necessary as in~\cite{BMTapp} to consider alternative hypotheses.

We assume that for every (sufficiently large) $p\in(1,\infty)$, there exists a Banach space $\cB$ containing constant functions, with norm $\|\,\|_{\cB}$,
and constants $\delta>0$, $\alpha_0\in(0,1)$ and $C>0$ such that
\begin{itemize}

\parskip = -2pt
\item[\textbf{(A)}]
\begin{itemize}
\item[(i)]
$\cB$ is compactly embedded in $L^p$.
\end{itemize}
\item[(ii)]
$\|\hat R(s)^n v\|_{\cB}\le C(|v|_{L^1}+\alpha_0^n \|v\|_{\cB})$
for all $s\in\barH_{\delta,L}$, $v\in\cB$, $n\ge1$.
\end{itemize}

It follows from these assumptions (see Lemma~\ref{lem:estKL}(d) below), that
(after possibly shrinking $\delta$) there is a continuous family of simple eigenvalues $\lambda(s)$ for $\hat R(s):\cB\to\cB$, $s\in\barH\cap B_\delta(0)$, with $\lambda(0)=1$.
Let $\zeta(s)\in\cB$ be the corresponding family of eigenfunctions normalized so that $\int_Y\zeta(s)\,d\mu=1$.
We assume further that there exists $\beta_+\in(\beta,1)$ such that
\begin{itemize}
\item[\textbf{(A)}]
(iii) $\big|\int_Y(e^{-s\tau}-1)(\zeta(s)-1)\,d\mu\big|\le C|s|^{\beta_+}$ for all $s\in \barH\cap B_\delta(0)$.
\end{itemize}

\begin{thm} \label{thm:alt}  
Suppose that hypothesis~(H) is replaced by hypothesis~(A).  Then
Theorems~\ref{thm:inf},~\ref{thm:average} and~\ref{thm:LLT} remain valid.  
If in addition
$\mu(\tau>t)=ct^{-\beta}+O(t^{-q})$ where $c>0$, $\beta\in(\frac12,1)$, $q>1$, 
then Theorem~\ref{thm:main} remains valid.
\end{thm}

\section{Mixing for infinite measure semiflows}
\label{sec:semiflow}

In this section, we obtain various mixing results for semiflows as consequences of the results in Section~\ref{sec:results}.

Let $F:Y\to Y$ and $\tau:Y\to\R^+$ be as in Section~\ref{sec:results}.
Define the suspension
$Y^\tau=\{(y,u)\in Y\times\R: 0\le u\le \tau(y)\}/\sim$
where $(y,\tau(y))\sim(Fy,0)$.  The suspension semiflow $F_t:Y^\tau\to Y^\tau$ is  given by $F_t(y,u)=(y,u+t)$, computed modulo identifications.
The measure
$\mu^\tau=\mu\times{\rm Lebesgue}$ is ergodic, $F_t$-invariant and $\sigma$-finite.  Since $\tau$ is nonintegrable, $\mu^\tau$ is an infinite measure.

Throughout this section, we suppose that
$A_1=A\times[a_1,a_2]$, $B_1=B\times[b_1,b_2]$ are measurable subsets of $\{(y,u)\in Y\times\R:0\le u\le \tau(u)\}$ 
(so $0\le a_1<a_2\le\essinf_A\tau$,
$0\le b_1<b_2\le\essinf_B\tau$),
and that $1_A\in\cB$.
Also, we continue to suppose that $\mu(\tau>t)=\ell(t)t^{-\beta}$ for various
ranges of $\beta\in[0,1]$.

\begin{cor}\label{cor:mixing}
Assume the setting of Theorem~\ref{thm:main} (alternatively Theorem~\ref{thm:alt}), so in particular $\beta\in(\frac12,1]$.  
Then
$\lim_{t\to\infty} m(t)\mu^\tau(A_1\cap F_t^{-1}B_1)  =d_\beta\mu^\tau(A_1)\mu^\tau(B_1)$.
\end{cor}

\begin{proof}
Recall that $\essinf\tau>1$.
Let $h\in(0,1)$ and note using~\eqref{eq:U} that
\begin{align*}
U_{A,B}(t+h)-U_{A,B}(t)
& = \mu(y\in A: F^ny\in B\;\text{and}\;\tau_n(y)\in[t,t+h]\; \text{for some}\; n\ge0)\\
&=\mu(y\in A: F_{t+h}(y,0)\in B\times[0,h]).
\end{align*}

After dividing rectangles into smaller subrectangles, we can suppose without loss that
$b_2-b_1<1$. Set $h=b_2-b_1$.  Then 
\begin{align}
\label{eq:rho}
\nonumber\mu^\tau(A_1\cap F_t^{-1}B_1) 
 & =\mu^\tau\{(y,u)\in A\times[a_1,a_2]:F_t(y,u)\in B\times[b_1,b_2]\}
 \\ \nonumber & =\mu^\tau\{(y,u)\in A\times[a_1,a_2]:F_{t+u-b_1}(y,0)\in B\times[0,h]\}
\\ \nonumber &  =\int_{a_1}^{a_2} \mu\{y\in A: F_{t+u-b_1}(y,0)\in B\times[0,h]\}\,du
 \\ & =\int_{a_1}^{a_2} (U_{A,B}(t+u-b_1)-U_{A,B}(t+u-b_1-h))\,du.
\end{align}
Hence
\begin{align*}
m(t)\mu^\tau(A_1\cap F_t^{-1}B_1) & = 
 \int_{a_1}^{a_2} m(t)(U_{A,B}(t+u-b_1)-U_{A,B}(t+u-b_1-h))\,du
 \\ & =\int_{a_1}^{a_2} \frac{m(t)}{m(t+u-b_1-h)}\,\chi(t+u-b_1-h)\,du,
\end{align*}
 where $\chi(t)=m(t)(U_{A,B}(t+h)-U_{A,B}(t))$ is bounded by Theorem~\ref{thm:main}.
Also $m(t)/m(t+u-b_1-h)$ is bounded by Potter's bounds (see for instance~\cite{BGT}).
Since $m(t)$ is regularly varying, we have $\lim_{t\to\infty}m(t)/m(t+u-b_1-h)=1$ for each $u\in[0,1]$.
By Theorem~\ref{thm:main}, $\lim_{t\to\infty}\chi(t+u-b_1-h)=d_\beta \mu(A)\mu(B) h= d_\beta\mu(A)\mu^\tau(B_1)$ for each $u\in[0,1]$.
Hence the result follows from the bounded convergence theorem.~
\end{proof}

\begin{rmk} \label{rmk:mixing}
The result also holds for all sets of the form $F_r^{-1}A_1$ and $F_s^{-1}B_1$
for fixed $r,s>0$.
Indeed, by Corollary~\ref{cor:mixing}, using that
$m(t)\sim m(t+s-r)$,
\begin{align*}
m(t)\mu^\tau(F_r^{-1}A_1\cap F_{t+s}^{-1}B_1)
  & = m(t)\mu^\tau(A_1\cap F_{t+s-r}^{-1}B_1)
\\ &  \to \mu^\tau(A_1)\mu^\tau(A_2)
= \mu^\tau(F_r^{-1}A_1)\mu^\tau(F_s^{-1}A_2).
\end{align*}
\end{rmk}

\begin{cor}[Weak rational ergodicity]
\label{cor:mixingav}
Assume the setting  of Theorem~\ref{thm:average} (alternatively Theorem~\ref{thm:alt}),
with $\beta\in[0,1]$.
Then
\[
\lim_{t\to\infty} t^{-1}m(t)\int_0^t \mu^\tau(A_1\cap F_x^{-1}B_1)\,dx =D_\beta\mu^\tau(A_1)\mu^\tau(B_1).
\]
\end{cor}

\begin{proof}
Continuing from~\eqref{eq:rho} (with $h=b_2-b_1$),
\begin{align*}
\int_0^t & \mu^\tau(A_1\cap F_x^{-1}B_1)\,dx  =\int_{a_1}^{a_2}\int_0^t (U_{A,B}(x+u-b_1)-U_{A,B}(x+u-b_1-h))\,dx\,du\\
&=\int_{a_1}^{a_2}\int_0^t U_{A,B}(x+u-b_1)\,dx\,du-\int_{a_1}^{a_2}\int_{-h}^{t-h} U_{A,B}(x+u-b_1)\,dx\,du\\
&=\int_{a_1}^{a_2}\int_{t-h}^t U_{A,B}(x+u-b_1)\,dx\,du-\int_{a_1}^{a_2}\int_{-h}^0 U_{A,B}(x+u-b_1)\,dx\,du
=I_1+I_2.
\end{align*}
Now
\[
t^{-1}m(t)I_1=t^{-1}m(t)U_{A,B}(t)\int_{a_1}^{a_2}\int_{-h}^0 \frac{U_{A,B}(x+t+u-b_1)}{U_{A,B}(t)}\,dx\,du.
\]
By Theorem~\ref{thm:average}, $U_{A,B}(t)$ is regularly varying so the integrand
$U_{A,B}(x+t+u-b_1)/U_{A,B}(t)$ is bounded for $x,u$ bounded and converges pointwise to $1$ as $t\to\infty$.  Hence
\[
\lim_{t\to\infty}\int_{a_1}^{a_2}\int_{-h}^0 \frac{U_{A,B}(x+t+u-b_1)}{U_{A,B}(t)}\,dx\,du=(a_2-a_1)h=(a_2-a_1)(b_2-b_1).
\]
By Theorem~\ref{thm:average}, $t^{-1}m(t)U_{A,B}(t)=D_{\beta}\mu(A)\mu(B)(1+o(1))$. 
Hence, $\lim_{t\to\infty}t^{-1}m(t)I_1=D_\beta
\mu(A)\mu(B)(a_2-a_1)(b_2-b_1)=\mu^\tau(A_1)\mu^\tau(B_1)$.
A simpler argument shows that $t^{-1}m(t)I_2=o(1)$. 
~\end{proof}

\begin{prop}  \label{prop:zerodensity}
Let $f:[0,\infty)\to\R$ be bounded and integrable on compact sets,
and let $K\in\R$.  Suppose that $\beta\in(0,1)$, that $\ell(t)$ is slowly
varying, and that
\begin{itemize}
\item[(a)] $\liminf_{t\to\infty}\ell(t)t^{1-\beta}f(t)\ge K$,
\item[(b)] $\lim_{t\to\infty}\ell(t)t^{-\beta}\int_0^t f(x)\, dx=\beta^{-1}K$.
\end{itemize}
Then there exists a set $E\subset[0,\infty)$ of density zero such that
$\lim_{t\to\infty,\;t\not\in E}\ell(t)t^{1-\beta}f(t)=K$.

In particular, 
$\liminf_{t\to\infty}\ell(t)t^{1-\beta}f(t)=K$.
\end{prop}

\begin{proof}
This is the continuous time analogue of~\cite[Proposition 8.2]{MT12} (which is itself a version of~\cite[p.~65, Lemma~6.2]{Petersen}).
We list the main steps which are proved exactly as in~\cite{MT12}.

\vspace{1ex}
\noindent{\bf Step~1}.
 Without loss of generality, $K=0$
and $\ell(t)t^{1-\beta}$ is increasing.

\vspace{1ex}
\noindent{\bf Step~2}.
Define the nested sequence of sets 
$E_q=\{t>0:\ell(t)t^{1-\beta}f(t)>1/q\}$, $q=1,2,\ldots$
Then $E_q$ has density zero for each $q$, i.e.\
$\lim_{t\to\infty}\frac{1}{t}\int_0^t 1_{E_q}(x)\, dx=0$.

\vspace{1ex}
\noindent{\bf Step~3}.
By Step 2, we can choose $0=i_0<i_1<i_2<\cdots$ such that
$\frac{1}{t}\int_0^t 1_{E_q}(x)\, dx<1/q$ for $t\ge i_{q-1}$, $q\ge2$.
Define $E=\bigcup_{q=1}^\infty E_q\cap(i_{q-1},i_q)$.
Then $E$ has density zero and 
$\lim_{t\to\infty,\;t\not\in E}\ell(t)t^{1-\beta}f(t)=0$.
\end{proof}

\begin{cor}\label{cor:inf}
Assume the setting of Theorem~\ref{thm:inf} (alternatively Theorem~\ref{thm:alt}),
with $\beta\in(0,1)$.  Then
\begin{itemize}
\item[(i)]
$\liminf_{t\to\infty} m(t)\mu^\tau(A_1\cap F_t^{-1}B_1)  =d_\beta\mu^\tau(A_1)\mu^\tau(B_1)$, 
and 
\item[(ii)] There exists a set $E\subset[0,\infty)$ of density zero such that \newline
$\lim_{t\to\infty,\;t\not\in E}m(t)\mu^\tau(A_1\cap F_t^{-1}B_1)  =d_\beta\mu^\tau(A_1)\mu^\tau(B_1)$.
\end{itemize}
\end{cor}

\begin{proof}
We start from the conclusion of Theorem~\ref{thm:inf}.
Arguing as in the proof of Corollary~\ref{cor:mixing}, but with $\lim$ replaced by $\liminf$ and using Fatou's lemma instead of the bounded convergence theorem, we obtain
\[
\liminf_{t\to\infty} \ell(t)t^{1-\beta}\mu^\tau(A_1\cap F_t^{-1}B_1)  \ge d_\beta\mu^\tau(A_1)\mu^\tau(B_1).
\]
This is condition~(a) in Proposition~\ref{prop:zerodensity},
and Corollary~\ref{cor:mixingav} 
is condition~(b).
Hence the result follows from
Proposition~\ref{prop:zerodensity}.
\end{proof}

\section{Main results used in the proof of Theorem~\ref{thm:main}}
\label{sec:steps}

The first result needed in  the proof of the strong renewal theorem, Theorem~\ref{thm:main},  is an inversion formula for the symmetric measure	
\[
\SMALL V_{A,B}(I)=\frac{1}{2}(U_{A,B}(I)+U_{A,B}(-I)).
\] 
Here, $U(-I)=U(\{x: -x\in I \})$
(with $U(-I)=0$ if $I\subset [0,\infty]$). We find it convenient
to adapt the formulation in~\cite[Section 4]{Erickson70}, but such an inversion
formula goes back to~\cite{FellerOrey61} (see also~\cite[Chapter 10]{Breiman}).\footnote{The result does not require any regular variation assumptions on $\mu(\tau>t)$, but we use the extra structure for simplicity.}

By (H) and (S)(i), $\hat T(s)=(I-\hat R(s))^{-1}$
is a bounded operator on $\cB$ for all $s\in\barH\setminus\{0\}$.
Let $A,\,B\subset Y$ be measurable with $1_A\in\cB$.

\begin{prop}[ {Analogue of~\cite[Inversion formula, Section 4]{Erickson70}.} ]
\label{prop:inv}
Let $g:\R\to\R$ be a continuous compactly supported function
with Fourier transform
$\hat g(x)=\int_{-\infty}^\infty e^{ixb} g(b)\, db$
satisfying
$\hat g(x)=O(x^{-2})$ as $x\to\infty$.
Then for all $\lambda,t\in\R$,
\[
\int_{-\infty}^\infty e^{-i\lambda(x-t)}\hat g(x-t)\,dV_{A,B}(x)=\int_{-\infty}^\infty e^{-itb} g(b+\lambda)\Re  \int_B \hat T(ib)1_A\,d\mu\, db.
\]
\end{prop}

The second result required in  the proof of Theorem~\ref{thm:main} comes directly from~\cite{Erickson70}
and does not require any modification in our set up. To state this result, for each $a>0$ we let $\hat g_a(0)=1$
and for $x\neq 0$, define
\[
\hat g_a(x)=\frac{2(1-\cos ax)}{a^2x^2}.
\]

\begin{prop}[ {\cite[Lemma 8]{Erickson70}} ]
\label{prop:Erickson}
Let $\{\mu_t,\,t>0\}$ be a family of measures such that $\mu_t(I)<\infty$ for every
compact set $I$ and all $t$. Suppose that for some constant $C$,
\[
\lim_{t\to\infty}\int_{-\infty}^{\infty} e^{-i\lambda x}\hat g_a(x)\, d\mu_t(x)=
C\int_{-\infty}^{\infty} e^{-i\lambda x}\hat g_a(x)\, dx,
\]
for all $a>0$, $\lambda\in\R$. Then $\mu_t(I)\to C |I|$
for every bounded interval $I$, where $|I|$ denotes the length of $I$. \qed
\end{prop}

Next, note that
$\hat g_a$ is the Fourier transform of
\[
g_a(b)=\begin{cases} a^{-1}(1-|b|/a), &  |b|\leq a \\ \hphantom{YYYY} 0, & |b|>a
\end{cases}.
\]
The final result required in  the proof of Theorem~\ref{thm:main} is as follows.

\begin{prop} \label{prop:Fga}
 For all $a>0$ and $\lambda\in\R$,
\[
\lim_{t\to\infty}m(t)\int_{-\infty}^\infty e^{-itb} g_a(b+\lambda)\Re  \int_B \hat T(ib) 1_A\,d\mu\, db=\pi d_\beta g_a(\lambda)\mu(A)\mu(B).
\]
\end{prop}

\begin{pfof}{Theorem~\ref{thm:main}}
With the convention $I+t=\{x:x-t\in I\}$, let  
\[
\mu_t(I)=2m(t)V_{A,B}(I+t)=m(t)(U_{A,B}(I+t)+U_{A,B}(-I-t))\]
and note that for $I=[0,h]$ with $h>0$,
\[
m(t)(U_{A,B}(t+h)-U_{A,B}(t))=\mu_t(I).
\]

Now,
\begin{align*}
m(t)\int_{-\infty}^\infty e^{-i\lambda(x-t)}\hat g_a(x-t)\,dV_{A,B}(x)&=
m(t)\int_{-\infty}^\infty e^{-i\lambda x}\hat g_a(x)\,dV_{A,B}(x+t)\\
&=\frac12\int_{-\infty}^{\infty} e^{-i\lambda x}\hat g_a(x)\, d\mu_t(x).
\end{align*}
Since $\hat g_a$ satisfies the assumptions of Proposition~\ref{prop:inv},
\[
\int_{-\infty}^{\infty} e^{-i\lambda x}\hat g_a(x)\, d\mu_t(x)=2m(t)\int_{-\infty}^\infty e^{-itb} g_a(b+\lambda)\Re  \int_B \hat T(ib) 1_A\,d\mu\, db.
\]
By Proposition~\ref{prop:Fga} together with the
Fourier inversion formula
$\int_{-\infty}^{\infty} e^{-i\lambda x}\hat g_a(x)\, dx=2\pi g_a(\lambda)$,
\[
\lim_{t\to\infty}\int_{-\infty}^{\infty} e^{-i\lambda x}\hat g_a(x)\, d\mu_t(x)
= 2\pi d_\beta g_a(\lambda)\mu(A)\mu(B) =d_\beta \int_{-\infty}^{\infty} e^{-i\lambda x}\hat g_a(x)\, dx\,\mu(A)\mu(B).
\]
Hence, we have shown that the hypothesis of Proposition~\ref{prop:Erickson} holds  with $C=d_\beta\mu(A)\mu(B)$.
It now follows from Proposition~\ref{prop:Erickson} with $I=[0,h]$ that
\[
m(t)(U_{A,B}(t+h)-U_{A,B}(t))=\mu_t([0,h])\to d_\beta \mu(A)\mu(B)h,
\]
as $t\to\infty$.
\end{pfof}

The proof of Propositions~\ref{prop:inv} and~\ref{prop:Fga} are
given in Section~\ref{sec:complete}.

\section{Prerequisites from operator renewal theory}
\label{sec:BMT}

In this section, we establish some
estimates for $\hat T =(I-\hat R)^{-1}$.
The arguments closely follow~\cite[Section~4]{BMTapp} (which was restricted to the case $\ell(t)=c+o(1)$ for some constant $c>0$ and did not include the case $\beta=1$). 

The estimates are carried out under hypotheses~(H) and~(S)(i) in Subsection~\ref{sec:H}.
The analogous results required under hypotheses~(A) and~(S)(i) are  obtained
in Subsection~\ref{sec:A}.

\subsection{Estimates under hypotheses~(H) and~(S)(i)}
\label{sec:H}

Throughout this subsection, $\beta\in(0,1]$ and $L>0$ are fixed.
We begin with $\gamma\in(0,\beta)$, $\delta\in(0,L)$ and 
$p>1$ as in (H).  
During the subsection, the values of $\gamma$, $\delta$ and $p$ change finitely many times; the changes in $\gamma$ are arbitrarily small.
Also $C>0$ is a constant whose value changes finitely many times.

For $r\in[0,1]$, let $\hat T_r(s)=(I-r\hat R(s))^{-1}$.
Define 
\[
\tilde\ell(t) =\begin{cases} \hphantom{XX} \ell(t) & \beta<1 \\
\int_1^t \ell(s)s^{-1}\,ds & \beta=1 \end{cases}, \qquad
c_\beta= \begin{cases}
i\int_0^\infty e^{-i\sigma}\sigma^{-\beta}\,d\sigma &  \beta<1 \\
\hphantom{XXX} 1 & \beta=1 \end{cases}.
\]

\begin{lemma}\label{lem:estKL}
(a) 
$\|\hat R(s_1)-\hat R(s_2)\|_{\cB\to L^p}\le C\, |s_1-s_2|^\gamma$
for all $s_1,s_2\in\barH$.
\\[.75ex]
(b)
There exists $r_0<1$ such that 
$\|\hat T_r(ib)\|_{\cB}\le C$
for all $|b|\in[\delta,L]$, $r\in[r_0,1]$.
\\[.75ex]
(c)
For all $\delta\le b<b+h<L$,
\[
\|\hat T(i(b+h))-\hat T(ib)\|_{\cB\to L^p}\le Ch^\gamma.
\]
(d) There exists a continuous family $\lambda(s)$, $s\in\barH\cap B_\delta(0)$, of simple eigenvalues for $\hat R(s):\cB\to\cB$ with $\lambda(0)=1$.
In addition, the corresponding family of spectral projections $P(s)$ are bounded linear operators on $\cB$ for all
$s\in\barH\cap B_\delta(0)$ and $\sup_{s\in\barH\cap B_\delta(0)}\|P(s)\|_{\cB}<\infty$.
Moreover, 
\[
\|P(s_1)-P(s_2)\|_{\cB\to L^p}
\le C |s_1-s_2|^\gamma
\quad\text{for all $s_1,s_2\in\barH\cap B_\delta(0)$.}
\]
(e) Define the complementary projections $Q(s)=I-P(s)$.  Then
 \[
\|(I-r\hat R(ib))^{-1}Q(ib)\|_\cB \le C
\quad \text{ for all $|b|<\delta$, $r\in[0,1]$.}
\]
\end{lemma}

\begin{proof}
(a)  Recall that $\hat R(s)v=R(e^{-s\tau}v)$.  Since $R$ is a positive operator,
\[
|(\hat R(s_1) -\hat R(s_2))v| 
 \le  R(|e^{-s_1\tau}-e^{-s_2\tau}||v|)
  \le  2|s_1-s_2|^\gamma R(\tau^\gamma |v|).
\]
By (H)(iii), $|(\hat R(s_1) -\hat R(s_2))v|_p\le 2|s_1-s_2|^\gamma
|R(\tau^\gamma |v|)|_p\ll |s_1-s_2|^\gamma \|v\|_{\cB}$.
\\[.75ex]
(b,c) 
Fix $b>0$.  It is immediate from
hypothesis~(S)(i) that $\|\hat T(ib)\|_\cB<\infty$.
Using also part~(a),
it follows from (H)(i,ii) and~\cite[Theorem~1]{KellerLiverani99}
that there exists $h_0>0$, $r_0<1$ 
and $C>0$ such that
$\|\hat T_r(i(b+h))\|_\cB\le C$ and
$\|\hat T(i(b+h))  -\hat T(ib)\|_{\cB\to L^p} \le C  |h|^\gamma$
for all $|h|<h_0$, $r\in(r_0,1]$.  
 The desired estimates follow from compactness of $[\delta,L]$.
\\[.75ex]
(d) This follows from (H)(i,ii) by~\cite[Corollary~1]{KellerLiverani99} exactly as in~\cite[Lemma~4.4]{BMTapp} (with $\beta-\eps$ replaced by $\gamma$).
\\[.75ex]
(e) By (H)(i,ii) and~\cite[Corollary~2]{KellerLiverani99}, for $\delta>0$ sufficiently small there exists $\rho\in(0,1)$ such that $\|(rR(ib)Q(ib))^n\|_{\cB}\le\|(R(ib)Q(ib))^n\|_{\cB} \le
C\rho^n$ for all $|b|<\delta$, $n\ge0$.  
\end{proof}

Let $\zeta(s)$ denote the corresponding family of
eigenfunctions normalized so that $\int_Y \zeta(s)\,d\mu=1$.
We have $\zeta(0)\equiv1$ and $P(0)v=\int_Y v\,d\mu$ for all $v\in \cB$.
It is immediate that $\zeta(s)$ 
inherits the estimates obtained for $P(s)$.  In particular, there is a constant $C>0$ such that $|\zeta(s)-\zeta(0)|_p\le C |s|^\gamma$
for all $s\in\barH\cap B_\delta(0)$.

Following~\cite{Gouezel10b} (see~\cite[Equation~(4.2)]{BMTapp}),
\begin{align} \label{eq:ev-Gou}
\SMALL \lambda(s) &  
\SMALL =\int_{Y}e^{-s\tau}\,d\mu +\chi(s)
 \quad\text{where}\quad
\chi(s)  =\int_Y (e^{-s\tau}-1)(\zeta(s)-\zeta(0))\,d\mu.
\end{align}

From now on, we fix $\delta\in(0,1)$ so that all conclusions of Lemma~\ref{lem:estKL} hold.

\begin{prop} \label{prop:IJM}
Write $s=a+ib\in\barH$.
\begin{itemize}
\item[(a)]
$1-\int_Y e^{-s\tau}\,d\mu\sim c_\beta \tilde\ell(1/|s|)s^\beta$ as $s\to0$.
\item[(b)]
When $\beta=1$,
$\Re\big(1-\int_Y e^{-ib\tau}\,d\mu\big)\sim \frac{\pi}{2}\ell(1/|b|)|b|$ as $b\to0$.
\item[(c)]  $|\int_Y (e^{-i(b+h)\tau}-e^{-ib\tau})d\mu|
\le C \tilde\ell(1/h)h^\beta$ for $0<h<b<\delta$.
\end{itemize}
\end{prop}

\begin{proof}
Part~(a) is proved as in~\cite[Lemma~2.4]{MT13} for $\beta<1$. Suppose that $\beta=1$ and 
let $G(x)=\mu(\tau>x)$.  Then
$1-\int_Y e^{-s\tau}\,d\mu = s\int_0^\infty e^{-sx}(1-G(x))\,dx
= sI_C(s)-isI_S(s)$, where
\[
\SMALL I_C(s)=\int_0^\infty e^{-ax}\cos bx\, (1-G(x))\,dx, \quad
I_S(s)=\int_0^\infty e^{-ax}\sin bx\, (1-G(x))\,dx.
\]

By~\cite[Proposition~6.2]{MT12}, we have for $a\ge |b|$ that
\[
I_C(s)=\tilde\ell(1/a)(1+o(1))+O(|b|a^{-1}\ell(1/a))
=\tilde\ell(1/|s|)(1+o(1))+O(\ell(1/|s|))\sim\tilde\ell(1/|s|).
\]
Similarly, for $a\le |b|$, we have
$I_C(s)=\tilde\ell(1/|b|)(1+o(1))+O(a|b|^{-1}\ell(1/|b|))
\sim\tilde\ell(1/|s|)$.
Hence $I_C(s)\sim \tilde\ell(1/|s|)$ as $s\to0$.
In the same way, it follows from~\cite[Proposition~6.2]{MT12} that $|I_S(s)|\ll \ell(1/|s|)$.
Part~(a) for $\beta=1$ follows immediately from these estimates.
Moreover, $I_S(ib)\sim\frac{\pi}{2}\ell(1/|b|)\sgn b$ as $b\to0$ by the
proof of~\cite[Lemma~6.8]{MT12}.
Since $\Re\big(1-\int_Y e^{-ib\tau}\,d\mu\big)=bI_S(ib)$, part~(b) follows.

Finally, part~(c) follows by the argument used in the proof of~\cite[Lemma 3.3.2]{GarsiaLamperti62}.~
\end{proof}

\begin{prop} \label{prop:chi}
\begin{itemize}
\item[(a)] $|\chi(s)|\le C|s|^{\beta+\gamma}$ for $s\in\barH\cap B_\delta(0)$,
\item[(b)] When $\beta>\frac12$, $|\chi(i(b+h))-\chi(ib)|\le C b^\beta h^\gamma$ for $0<h<b<\delta$.
\end{itemize}
\end{prop}

\begin{proof}
Choose $\eps>0$ arbitrarily small and $r>1$ such that $(\beta-\eps)r<\beta$ with conjugate exponent $r'$.
Then $\tau^{(\beta-\eps)r}\in L^1$ and it follows from
H\"older's inequality that
\[
\SMALL |\chi(s)| \le 2|s|^{\beta-\eps}| \tau^{\beta-\eps}(\zeta(s)-1)|_1
\le 2|s|^{\beta-\eps}|\tau^{\beta-\eps}|_r|\zeta(s)-1|_{r'}
\ll |s|^{\beta-\eps+\gamma},
\]
yielding part~(a).   Here we used that $|\zeta(s)-1|_p=O(|s|^\gamma)$ for $p$ as large as desired.
Similarly,
\begin{align*}
|\chi(i(b+h))-\chi(ib)|
& \le |(e^{i(b+h)\tau}-1)(\zeta(i(b+h))-\zeta(ib))|_1
+ |(e^{ih\tau}-1)(\zeta(ib)-1)|_1
\\ & \ll (b+h)^{\beta-\eps}h^\gamma+h^{\beta-\eps}b^\gamma
 \ll b^\beta h^{\gamma-\eps},
\end{align*}
(Note that $h^{\beta-\eps}b^\gamma=h^{\gamma-\eps}h^{\beta-\gamma}b^\gamma
\le h^{\gamma-\eps}b^\beta$ since $\gamma<\beta$ and $h<b$.)
This proves part~(b).
\end{proof}

\begin{cor} \label{cor:lambda}
Write $s=a+ib\in\barH$.
\begin{itemize}
\item[(a)]
$1-\lambda(s)\sim  c_\beta\tilde\ell(1/|s|)s^{\beta}$ as $s\to0$.
\item[(b)]
When $\beta=1$, 
$\Re(1-\lambda(ib))\sim  \frac{\pi}{2}\ell(1/|b|)|b|$ as $b\to0$. 
\item[(c)]
When $\beta>\frac12$,
$|\lambda(i(b+h))-\lambda(ib)| \le C(
\tilde\ell(1/h)h^\beta +b^\beta h^\gamma)$,
 for $0<h< b<\delta$.
\item[(d)] $|1-r\lambda(ib)|^{-1}\le C\tilde\ell(1/|b|)^{-1}|b|^{-\beta}$ for all $|b|<\delta$, $r\in[\frac12,1]$.
\item[(e)] When $\beta=1$, $|\Re(1-r\lambda(ib))|^{-1}\le C\ell(1/|b|)\tilde\ell(1/|b|)^{-2}|b|^{-1}$ for all $|b|<\delta$, $r\in[\frac12,1]$.
\end{itemize}
\end{cor}

\begin{proof}
Parts~(a) and (b) are immediate from~\eqref{eq:ev-Gou} and Propositions~\ref{prop:IJM}(a,b) and~\ref{prop:chi}(a).
Part~(c) follows from~\eqref{eq:ev-Gou} and Propositions~\ref{prop:IJM}(c) and~\ref{prop:chi}(b).
By part~(a),
\[
\SMALL |1-r\lambda(ib)|\ge |\Im(1-r\lambda(ib))|\ge \frac12 |\Im\lambda(ib)|\sim
\frac12 |\Im(i^\beta c_\beta)|\tilde\ell(1/|b|)|b|^\beta,
\]
yielding part~(d).
Using also that $\Re\lambda(ib)\in[0,1]$ for $|b|<\delta$, we compute for
$\beta=1$ that
\begin{align*}
|\Re(1-r\lambda(ib)^{-1}|
& =\Re(1-r\lambda(ib))\,|1-r\lambda(ib)|^{-2}
\\ & \le \Re(1-\lambda(ib))\,|r\Im(\lambda(ib))|^{-2}
 \le 4\Re(1-\lambda(ib))\,|r\Im(\lambda(ib))|^{-2},
\end{align*}
so part~(e) follows from parts (a) and (b).
\end{proof}

\begin{lemma} \label{lem:T}
  $\hat T(s)  =c_\beta^{-1} \tilde\ell(1/|s|)^{-1}s^{-\beta}(P(0) +E(s))$
for $s\in\barH\cap B_\delta(0)$, 
where  $E(s)$ is a family of operators satisfying $\lim_{s\to0}\|E(s)\|_{\cB\to L^1}=0$.
\end{lemma}

\begin{proof}
By Corollary~\ref{cor:lambda}(a),
$(1-\lambda(s))^{-1}\sim  c_\beta^{-1}\tilde\ell(1/|s|)^{-1}s^{-\beta}$ as $s\to0$.
Also,
\[
\hat T(s)  =(1-\lambda(s))^{-1}P(s)+(I-\hat R(s))^{-1}Q(s)
 =(1-\lambda(s))^{-1}(P(0)+E(s)),
\]
where
\begin{align} \label{eq:E}
 E(s)= P(s)-P(0)+(1-\lambda(s))(I-\hat R(s))^{-1}Q(s).
\end{align}
By (H), $\|(I-\hat R(s))^{-1}Q(s)\|_{\cB}=O(1)$.
By  Lemma~\ref{lem:estKL}(d), $\|P(s)-P(0)\|_{\cB\to L^1}=O(|s|^\gamma)$.  Hence
$\|E(s)\|_{\cB\to L^1}\ll |s|^\gamma+|s|^{\beta-\eps}$.
\end{proof}

\begin{lemma} \label{lem:1}
Let $\beta=1$.  Then
  $\Re \hat T(ib)  = \frac{\pi}{2} \ell(1/|b|)\tilde\ell(1/|b|)^{-2}|b|^{-1}(P(0) +E(b))$
for $b\in\R$, $0<|b|<\delta$,
where  $\lim_{b\to0}\|E(b)\|_{\cB\to L^1}=0$.
\end{lemma}

\begin{proof}
By Corollary~\ref{cor:lambda}(a,b),
\[
\Re ((1-\lambda(ib))^{-1})= \Re(1-\lambda(ib))|1-\lambda(ib)|^{-2}
\sim \frac{\pi}{2}\ell(1/|b|)\tilde\ell(1/|b|)^{-2}|b|^{-1}.
\]
As in the proof of Lemma~\ref{lem:T}, $\Re\hat T(ib)=
\{\Re((1-\lambda(ib))^{-1})\}(P(0)+E(b))$ where
\[
E(b)=\Re(1-\lambda(ib))\Re\{(1-\lambda(ib))^{-1}(P(ib)-P(0))+(I-R(ib))^{-1}Q(ib)\},
\]
and
\[
\|E(b)\|_{\cB\to L^1}  \ll 
\|P(ib)-P(0)\|_{\cB\to L^1}+
|1-\lambda(ib)|\|(I-R(ib))^{-1}Q(ib)\|_{\cB\to L^1})\ll |b|^{1-\eps},
\]
completing the proof.
\end{proof}

\begin{cor} \label{cor:psi}
Let $\beta\le1$, $L>0$.  There are constants $r_0<1$ and $C>0$ such that
\[
\|\Re\hat T_r(ib)\|_{\cB\to L^1}\le C \psi_\beta(|b|)
\quad\text{for $0<|b|\le L$, $r_0\le r\le 1$,} 
\]
where
  \[
\psi_\beta(x)= \begin{cases}  \ell(1/x)^{-1}x^{-\beta} & \beta<1 \\ \ell(1/x)\tilde\ell(1/x)^{-2}x^{-1} & \beta=1 \end{cases}.
\]
\end{cor}

\begin{proof}
By Lemma~\ref{lem:estKL}(b), we can restrict to the range $|b|<\delta$ on which
\[
\hat T_r(ib)=(I-r\hat R(ib))^{-1}=(1-r\lambda(ib))^{-1}P(ib)+(I-r\hat R(ib))^{-1}Q(ib).
\]
The result follows from the estimates for $P$, $(I-r\hat R)^{-1}Q$ and 
$(1-r\lambda)^{-1}$ obtained in Lemma~\ref{lem:estKL}(d,e) and
Corollary~\ref{cor:lambda}(d,e).
\end{proof}

\begin{rmk} \label{rmk:psi}
(a) 
Note that $\psi_\beta$ is integrable on $[0,L]$ for all $\beta\le1$.
This is clear for $\beta<1$ while $\tilde\ell(1/x)^{-1}$ is an antiderivative for
$\psi_1$.
In particular, $\sup_r|\Re \hat T_r(ib)1_A|_1\le C \psi_\beta(b)\|1_A\|_\cB$ which is integrable.
\\[.75ex]
(b) By Karamata's theorem on integration of regularly varying sequences~\cite{BGT}, $\tilde\ell$ is slowly varying and
$\ell(x)=o(\tilde\ell(x))$ as $x\to\infty$ when $\beta=1$.  In particular,
$\psi_\beta(b)\ll \tilde\ell(1/|b|)^{-1}|b|^{-\beta}$ for all $\beta\le1$.
\end{rmk}

\begin{lemma}
\label{lem:Tcont}
Let $\beta\in(\frac12,1]$.
 For $0<h< b<\delta$,
  \[
 \|\hat T(i(b+h))  -\hat T(ib)\|_{\cB\to L^1}
\le C\big\{ \tilde\ell(1/b)^{-2}b^{-2\beta}\tilde\ell(1/h) h^\beta + b^{-\beta} h^\gamma \big\}.
\]
\end{lemma}

\begin{proof}
Recall as in Lemma~\ref{lem:T} that
$\hat T(ib)  =A_1(b)+A_2(b)$, where
\[
A_1(b)  =(1-\lambda(ib))^{-1}P(ib), \qquad A_2(b)=(I-\hat R(ib))^{-1}Q(ib).
\]
Using Lemma~\ref{lem:estKL}(d) and Corollary~\ref{cor:lambda}(a,c),
\begin{align*}
\|A_1(b+h)-A_1(b)\|_{\cB\to L^1} & \ll 
|1-\lambda(i(b+h))|^{-1}\|P(i(b+h))-P(ib)\|_{\cB\to L^1}
\\ & \!\!\!\!\!\!\!\! 
 \!\!\!\!\!\!\!\! 
 \!\!\!\!\!\!\!\! 
 \!\!\!\!\!\!\!\! 
+|1-\lambda(ib)|^{-1}|1-\lambda(i(b+h))|^{-1}|\lambda(i(b+h))-\lambda(ib)|\|P(ib)\|_{\cB\to L^1}
\\ & \ll
\tilde\ell(1/b) b^{-\beta}h^{\gamma}
+\tilde\ell(1/b)^{-2}b^{-2\beta}(\tilde\ell(1/h)h^\beta 
+b^\beta h^\gamma)
\\ & \ll \tilde\ell(1/b)^{-2}b^{-2\beta}\tilde\ell(1/h)h^\beta + b^{-\beta}h^{\gamma-\eps}.
\end{align*}
An argument
from~\cite[Proposition 3.8]{LiveraniTerhesiu16} shows that
$\|A_2(b+h)-A_2(b)\|_{\cB\to L^1}\ll h^{\gamma-\eps}$, completing the proof.
\end{proof}

\subsection{Estimates under hypotheses~(A) and~(S)(i)}
\label{sec:A}

Let $\eps\in(0,\beta)$.
Since $R:L^1\to L^1$ is a contraction,
\[
|(\hat R(s_1) -\hat R(s_2))v|_1 
 \le  |(e^{-s_1\tau}-e^{-s_2\tau})v|_1
  \le  2|s_1-s_2|^{\beta-\eps}|\tau^{\beta-\eps} v|_1.
\]
Choose $r>1$ such that $(\beta-\eps)r<\beta$ with conjugate exponent $r'$.
By H\"older's inequality and (A)(i),
$|\tau^{\beta-\eps} v|_1\le |\tau^{\beta-\eps}|_r|v|_{r'}\ll \|v\|_{\cB}$.
Hence
$\|\hat R(s_1)-\hat R(s_2)\|_{\cB\to L^1}\ll |s_1-s_2|^{\beta-\eps}$
for all $s_1,s_2\in\barH$.  

Using~\cite{KellerLiverani99} as before, we deduce that
the conclusions of Lemma~\ref{lem:estKL} hold with $L^p$ replaced by $L^1$
and $\gamma$ replaced by $\beta-\eps$.

\begin{prop} \label{prop:psi}
The conclusions of Lemmas~\ref{lem:T} and~\ref{lem:1} and 
Corollary~\ref{cor:psi} are unchanged under hypotheses~(A) and (S)(i).
\end{prop}

\begin{proof}
It is immediate from hypothesis~(A)(iii) that $|\chi(s)|\ll |s|^{\beta_+}$ where $\beta_+>\beta$, and hence the proofs 
are unchanged.
\end{proof}

Lemma~\ref{lem:Tcont} becomes:

\begin{lemma} \label{lem:TcontA}
 $\|\hat T(i(b+h)) -\hat T(ib)\|_{\cB\to L^1}
\le Cb^{-2\beta} h^{\beta-\eps}$ for all  $0<h< b<\delta$.
\end{lemma}

\begin{proof}
Since $\|\zeta(s)\|_{\cB}$ is bounded, it follows again from H\"older's inequality that
\begin{align*}
|\chi(i(b+h))-\chi(ib)|
& \le |(e^{i(b+h)\tau}-1)(\zeta(i(b+h))-\zeta(ib))|_1
+|(e^{ih\tau}-1)(\zeta(ib)-1)|_1  \\
& \le 2|\zeta(i(b+h))-\zeta(ib)|_1+
2h^{\beta-\eps} |\tau^{\beta-\eps}|_r|(\zeta(ib)-1)|_{r'}
 \ll h^{\beta-\eps}.
\end{align*}
Hence by~\eqref{eq:ev-Gou} and Proposition~\ref{prop:IJM}(c),
$|\lambda(i(b+h))-\lambda(ib)| \ll h^{\beta-\eps}$.
Now proceed as in the proof of Lemma~\ref{lem:Tcont}.
\end{proof}

The presence of the $\eps$ in Lemma~\ref{lem:TcontA} necessitates some alterations to the strategy in~\cite{Erickson70}.   As in~\cite{BMTapp}, we make use of the following refinement of Lemma~\ref{lem:T}.

\begin{lemma}
\label{lem:T-higherorder}
Assume that $\mu(\tau>t)=ct^{-\beta}+O(t^{-q})$ where $c>0$, $\beta\in(\frac12,1)$, $q>1$.  
Then
  $c\hat T(ib)  =c_\beta^{-1}  b^{-\beta}P(0) +\tilde E(b)$
for $b\in [0,\delta)$, 
where
$\|\tilde E(b)\|_{\cB\to L^1}\le Cb^{-(2\beta-\beta_+)}$.
\end{lemma}

\begin{proof}
A calculation using only the expression for $\mu(\tau>t)$ shows that 
$1-\int_Y e^{-s\tau}\,d\mu=cc_\beta b^\beta +O(b)$  
(see~\cite[Eq.\ (4.4)]{BMTapp}).
By~\eqref{eq:ev-Gou} and the estimate $|\chi(s)|\ll |s|^{\beta_+}$ where
$\beta_+\in(\beta,1)$, we obtain 
$1-\lambda(ib)=cc_\beta b^\beta(1+O(b^{\beta_+-\beta}))$.
Hence
\[
c(1-\lambda(ib))^{-1}=c_\beta^{-1} b^{-\beta}(1+O(b^{-(2\beta-\beta_+)})).
\]
By~\eqref{eq:E}, 
$c\hat T(ib)=c(1-\lambda(ib))^{-1}P(0)+c(1-\lambda(ib))^{-1}E(ib)
=c_\beta^{-1} b^{-\beta}P(0)+c\tilde E(b)$ where
$\tilde E(b)=(1-\lambda(ib))^{-1}E(ib)+O(b^{-(2\beta-\beta_+)})$ and
\[
E(ib)= P(ib)-P(0)+(1-\lambda(ib))(I-\hat R(ib))^{-1}Q(ib)=O(b^{\beta-\eps}).
\]
Hence $\tilde E(b)\ll b^{-\eps}+b^{-(2\beta-\beta_+)}$.  
Recall that $2\beta-\beta_+>2\beta-1>0$, so we can
choose $\eps\in(0,2\beta-\beta_+)$ completing the proof.
\end{proof}

\section{Completion of the proof of Theorem~\ref{thm:main}}
\label{sec:complete}

In this section, we give the proof of Propositions~\ref{prop:inv} and~\ref{prop:Fga}, thereby completing the proof of Theorem~\ref{thm:main}.
In Subsections~\ref{sec:inv} and~\ref{sec:Fga}, we assume hypotheses~(H) and~(S)(i).
In Subsection~\ref{sec:completeA}, we show that the results remain true under hypotheses~(A) and~(S)(i).

\subsection{Proof of Proposition~\ref{prop:inv}}
\label{sec:inv}

Fix $\beta\le1$.
Throughout, we write $U$ and $V$ instead of $U_{A,B}$ and $V_{A,B}$.
Following~\cite[Chapter~10]{Breiman} (see also~\cite[Section~4]{Erickson70}),
we define for $r\in(0,1)$,
\begin{align*}
U_r(I) & = \sum_{n=0}^\infty r^n\mu(y\in A\cap F^{-n}B:\tau_n(y)\in I),
\\
V_r(I) & = \SMALL\frac{1}{2}(U_r(I)+U_r(-I)).
\end{align*}

For $n\ge0$, the Fourier transform of the distribution $G_n(x)=\mu(\tau_n(y)\leq x,\,y\in A\cap F^{-n}B)$
is given by $\int_Y 1_A\,1_B\circ F^n\,e^{ib\tau_n}\, d\mu=\int_B \hat R(-ib)^n 1_A\,d\mu$.
 Hence
\begin{align*}
\nonumber \int_{-\infty}^{\infty} e^{ibx}\, d V_r(x)&=\Re\int_{0}^{\infty} e^{ibx}\, d U_r(x)
\\ & =\sum_{n=0}^\infty r^n\Re\int_B \hat R(ib)^n 1_A\,d\mu
=\Re\int_B \hat T_r(ib) 1_A\,d\mu,
\end{align*}
where $\hat T_r(s)=(I-r\hat R(s))^{-1}$.

Let $\hat g$ and $g$ be as in the statement of Proposition~\ref{prop:inv}. 
Note that $dV_r$ is a finite measure and $g$ is compactly supported, so $e^{ibx}g(b)$ lies in $L^1(dV_r\times db)$.  Hence
it follows from Fubini's theorem that for $r\in(0,1)$,
\begin{align*}
 \int_{-\infty}^{\infty} &  \hat g(x)\, d V_r(x)  =  
 \int_{-\infty}^{\infty}\Big\{ \int_{-\infty}^{\infty} e^{ibx}g(b)\,db\Big\}\, d V_r(x)
\\ & =  
 \int_{-\infty}^{\infty} g(b)\Big\{ \int_{-\infty}^{\infty} e^{ibx}\,dV_r(x)\Big\}\, db
=
\int_{-\infty}^{\infty}g(b)\Re\int_B \hat T_r(ib) 1_A\,d\mu\, db.
\end{align*}

Replacing $g(b)$ by $g_1(b)=e^{-ibt}g(b+\lambda)$ and
$\hat g(x)$ by 
$\hat g_1(x)=\int_{-\infty}^{\infty}e^{ibx}\, g_1(b)\, db=e^{-i\lambda(x-t)}\hat g(x-t)$,
we obtain
\begin{align} \label{eq:Vr}
 \int_{-\infty}^{\infty} &  e^{-i\lambda(x-t)}\hat g(x-t)\, d V_r(x)  =  
\int_{-\infty}^{\infty}e^{-ibt}g(b+\lambda)\Re\int_B \hat T_r(ib) 1_A\,d\mu\, db.
\end{align}
It remains to justify passing to the limit $r\to1_-$ on both sides of~\eqref{eq:Vr}.

First, we consider the left-hand side of~\eqref{eq:Vr}.
Since $\tau\ge1$, we have $U(x)=U([0,x])\le 
\sum_{n=0}^\infty \mu(\tau_n\le x)\le x+1$ for all $x$.  Integrating by parts,
\[
\int_1^\infty x^{-2}\,dU(x)=
-U(1)+2\int_1^\infty U(x)x^{-3}\,dx<\infty.
\]
Hence $\int_{|x|\ge1}x^{-2}dV(x)<\infty$.
Since $\int_{-1}^1|\hat g(x-t)|\,dV(x)<\infty$ and $\hat g(x-t)=O(x^{-2})$ for each fixed $t$, it follows that
$f(x)=e^{-i\lambda(x-t)}\hat g(x-t)$ is integrable with respect to $dV(x)$. 
But $V_r(I)\nearrow V(I)$ as $r\to1_-$ for every measurable $I$,
so $\lim_{r\to1_-}\int_{-\infty}^\infty f(x)\,dV_r(x)=
\int_{-\infty}^\infty f(x)\,dV(x)$ which is the required result for the left-hand side.

Finally, we consider the right-hand side of~\eqref{eq:Vr}.
Choose $L>0$ such that $\supp g\in[-L,L]$.
By Remark~\ref{rmk:psi}(a),
$|\Re\hat T_r(ib)1_A|_1\ll \psi_\beta(b)\|1_A\|_{\cB}$ for $|b|\le L+|\lambda|$, where
$\psi_\beta$ is integrable. 
Hence the desired limit as $r\to1_-$ follows from the dominated convergence theorem.

\subsection{Proof of Proposition~\ref{prop:Fga}}
\label{sec:Fga}

Fix $\beta\in(\frac12,1]$.
We follow the proof of~\cite[Theorem 1]{Erickson70} (an adaptation of the argument in~\cite{GarsiaLamperti62}).
Let $W(b)=\Re  \int_B \hat T(ib)1_A\,d\mu$.

Fix $\omega>1$ and write 
$\int_{-\infty}^\infty e^{-itb} g_a(b+\lambda)\Re  \int_B \hat T(ib) 1_A\,d\mu\, db=I_1(t,\omega)+I_2(t,\omega)$ where
\[
I_1(t,\omega)  = \int_{-\omega/t}^{\omega/t} e^{-itb}g_a(b+\lambda)W(b)\, db,  \quad
I_2(t,\omega)  = 
\int_{|b|>\omega/t} e^{-itb}g_a(b+\lambda)W(b)\, db.
\]

Proposition~\ref{prop:Fga} follows immediately from the estimates
for $I_1(t,\omega)$ and $I_2(t,\omega)$ below.

\begin{lemma} \label{lem:I1}
$\lim_{\omega\to\infty}\lim_{t\to\infty} m(t)I_1(t,\omega) = \pi d_\beta g_a(\lambda)\mu(A)\mu(B) $.
\end{lemma}

\begin{proof}
It follows from the definition of $g_a$ that $|g_a(b_1)-g_a(b_2)|\le a^{-2}|b_1-b_2|$.
Hence
\begin{align*}
\Big|I_1(t,\omega)-g_a(\lambda) \int_{-\omega/t}^{\omega/t} e^{-itb}W(b)\,db\Big| & \le \int_{-\omega/t}^{\omega/t} |g_a(b+\lambda)-g_a(\lambda)|\,|W(b)|\,db \\
& \le 2a^{-2}\omega t^{-1}\int_0^{\omega/t}|W(b)|\,db.
\end{align*}
By Remark~\ref{rmk:psi}(a), 
 $\int_0^{\omega/t}|W(b)|\,db \ll \|1_A\|$
for $t>\omega/\delta$.  Hence
\[
\lim_{t\to\infty}m(t)I_1(t,\omega)
 =2g_a(\lambda)\lim_{t\to\infty}m(t) \int_0^{\omega/t} W(b)\cos tb\,db.
\]

For $\beta<1$, define $\xi(b)=\mu(A)\mu(B)+\int_B E(ib)1_A\,d\mu$ where
$E$ is as in Lemma~\ref{lem:T}.
In particular,
$|\xi(b)|\le |1_A|_1+|E(ib)1_A|_1\ll \|1_A\|$ and 
$|\xi(b)-\mu(A)\mu(B)|\le \|E(ib)\|_{\cB\to L^1}\|1_A\|\to0$ as $b\to0$.
Hence
\begin{align*}
 m(t) \int_0^{\omega/t} W(b)\cos tb\,db
  & =\ell(t) t^{1-\beta}\Re\Big\{ c_\beta^{-1}\int_0^{\omega/t} \ell(1/b)^{-1}b^{-\beta}\xi(b)\cos tb \,db\Big\}
 \\ & 
 = \Re\Big\{ c_\beta^{-1}\int_0^{\omega} [\ell(t)/\ell(t/b)]b^{-\beta}\xi(b/t)\cos b \,db\Big\}.
\end{align*}
By the dominated convergence theorem,
\[
\lim_{t\to\infty}m(t)\int_0^{\omega/t} W(b)\cos tb\,db = (\Re c_\beta^{-1})\int_0^\omega b^{-\beta}\cos b\,db\,\mu(A)\mu(B),
\]
and the result for $\beta<1$ follows.  

Now suppose that $\beta=1$ and recall that $\psi_1(b)=\ell(1/b)\tilde\ell(1/b)^{-2}b^{-1}$.
By Lemma~\ref{lem:1},
\[
m(t) \int_0^{\omega/t} W(b)\cos tb\,db
= \tilde\ell(t)\frac{\pi}{2}\int_0^{\omega/t}\psi_1(b)\xi(b)\cos tb\,db.
\]
where $\xi(b)$ has the same properties as before.
Now
\begin{align*}
\tilde\ell(t)\int_0^{\omega/t}\psi_1(b)\xi(b)\,db
& =\tilde\ell(t)\int_0^{\omega/t}\psi_1(b)(\mu(A)\mu(B)+o(1))\,db
\\ & =\tilde\ell(t)\tilde\ell(t/\omega)^{-1}(\mu(A)\mu(B)+o(1))\to \mu(A)\mu(B).
\end{align*}
Next,
\begin{align*}
\tilde\ell(t)\int_0^{\omega/t}\psi_1(b)\xi(b)(\cos tb-1) \,db & =\int_0^\omega 
\frac{\tilde\ell(t)}{\tilde\ell(t/\sigma)}
\frac{\ell(t/\sigma)}{\tilde\ell(t/\sigma)} \xi(\sigma/t)
\frac{\cos \sigma-1}{\sigma}
\,d\sigma.
\end{align*}
By Remark~\ref{rmk:psi}(b), $\tilde\ell$ is slowly varying and $\ell(x)=o(\tilde\ell(x))$ as $x\to\infty$.
By Potter's bounds, the integrand is dominated by $\sigma^{1-\eps}$ for any $\eps>0$, so
the integrand converges to zero pointwise and
$\tilde\ell(t)\int_0^{\omega/t}\psi_1(b)\xi(b)(\cos tb-1)\,db \to0$ as $t\to\infty$.
Hence $\lim_{t\to\infty}m(t)\int_0^{\omega/t}W(b)\cos tb\,db=\frac{\pi}{2}\mu(A)\mu(B)$ yielding the result for $\beta=1$.
\end{proof}

\begin{lemma} \label{lem:I2}  Let $\beta'\in(\frac12,\beta)$.
Then $\limsup_{t\to\infty} m(t)I_2(t,\omega) = O(\omega^{-(2\beta'-1)})$.
\end{lemma}

\begin{proof}
It follows from evenness of $g_a$ and $W(b)$, together with the fact that $\supp g_a=[-a,a]$, that
\[
I_2(t,\omega)  = \int_{b>\omega/t} [e^{-itb}g_a(b+\lambda)+e^{itb}g_a(b-\lambda)]W(b)\, db
 = \int_{\omega/t}^{a+|\lambda|} h(b)W(b)\, db,
\]
where 
$h(b)=e^{-itb}g_a(b+\lambda)+e^{itb}g_a(b-\lambda)$.
Continuing as on~\cite[p.~278]{Erickson70} down as far as~\cite[Equation~(5.14)]{Erickson70}, we obtain
$m(t)|I_2(t,\omega)|  \le a^{-1}J_1(t,\omega)+ \pi a^{-2}J_2(t,\omega)+ a^{-1}J_3(t,\omega)$,
where
\begin{align*}
J_1(t,\omega) &= m(t)\int_{(\omega-\pi)/t}^{\omega/t}|W(b+\pi/t)|\,db , \quad
J_2(t,\omega) =m(t)t^{-1}\int_{\omega/t}^{a+|\lambda|}|W(b)|\,db , \\
J_3(t,\omega) & = m(t)\int_{\omega/t}^{a+|\lambda|}|W(b+\pi/t)-W(b)|\,db.
\end{align*}

By Remark~\ref{rmk:psi}(a), $W$ is integrable on $[0,a+|\lambda|]$ so 
$J_2(t,\omega)\ll \tilde\ell(t)t^{-\beta}\to0$ as $t\to\infty$.
By Lemma~\ref{lem:T},
for $\beta<1$,
\[
J_1(t,\omega)\ll \ell(t)t^{1-\beta}\int_{\omega/t}^{(\omega+\pi)/t}\ell(1/b)^{-1}b^{-\beta}\,db
= \int_{\omega}^{\omega+\pi}(\ell(t)/\ell(t/\sigma))\sigma^{-\beta}\,d\sigma\ll
\omega^{-(\beta-\eps)},
\]
for any $\eps>0$ by Potter's bounds.
By Lemma~\ref{lem:1} and Remark~\ref{rmk:psi}(b),
for $\beta=1$,
\[
J_1(t,\omega)\ll \tilde\ell(t)\int_{\omega/t}^{(\omega+\pi)/t} \psi_1(b)\,db=
\tilde\ell(t)\{\tilde\ell(t/(\omega+\pi))^{-1}-\tilde\ell(t/\omega)\}\to0
\quad\text{as $t\to\infty$}.
\]

By Lemma~\ref{lem:Tcont} with $h=\pi/t$,
\[
J_3(t,\omega)  \ll 
\tilde\ell(t)^2t^{1-2\beta}\int_{\omega/t}^\infty \tilde\ell(1/b)^{-2}b^{-2\beta}\,db+
t^{1-\beta+\eps-\gamma}\int_0^{a+|\lambda|}b^{-\beta}\,db 
=J_{3,1}+J_{3,2}.
\]
By Potter's bounds,
\[
J_{3,1}  
=\int_\omega^\infty  [\tilde\ell(t)/\tilde\ell (t/\sigma)]^2  \sigma^{-2\beta}\,d\sigma 
\ll  \int_\omega^\infty \sigma^{-2\beta'}\,d\sigma  
\ll
\omega^{-(2\beta'-1)}.
\]
Finally, since we are in the case $\beta>\frac12$, we can choose $\gamma\in(1-\beta,\beta)$ in hypothesis~(H).
Hence $J_{3,2} \ll t^{1-\beta+\eps-\gamma}=o(1)$ as $t\to\infty$
for $\eps>0$ sufficiently small.
\end{proof}

\subsection{Modified argument under hypotheses~(A) and~(S)(i)}
\label{sec:completeA}

Assume hypotheses~(A) and~(S)(i) and that
$\mu(\tau>t)=ct^{-\beta}+O(t^{-q})$ where $c>0$, $\beta\in(\frac12,1)$, $q>1$.
Recall that $\beta_+>\beta$.

First, we note by Proposition~\ref{prop:psi} that Corollary~\ref{cor:psi}
is unchanged.  Hence the proof of Proposition~\ref{prop:inv}
is unchanged.

For Proposition~\ref{prop:Fga}, we adopt a different strategy from before.
Instead of considering $\lim_{\omega\to\infty}\limsup_{t\to\infty}I_r(t,\omega)$ for $r=1,2$, we consider 
$\lim_{t\to\infty} I_r(t,t^{\kappa})$ for a suitable choice of $\kappa>0$.

\begin{lemma} \label{lem:I1A}
$\lim_{t\to\infty}m(t)I_1(t,t^{\kappa}) = \pi d_\beta g_a(\lambda)\mu(A)\mu(B)$
for all $\kappa>0$.
\end{lemma}

\begin{proof}
Following the proof of Lemma~\ref{lem:I1} and using Lemma~\ref{lem:T} and Proposition~\ref{prop:psi},
\[
\Big|m(t)I_1(t,\omega)-2m(t)g_a(\lambda)\int_0^{\omega/t}W(b)\cos tb\,db\Big|\ll \omega t^{-\beta}\int_0^{\omega/t}b^{-\beta}\,db \ll \omega^{2-\beta}t^{-1}.
\]
By Lemma~\ref{lem:T-higherorder},
\begin{align*}
m(t)\int_0^{\omega/t}W(b)\cos tb\,db 
& = t^{1-\beta}\int_0^{\omega/t} (\Re c_\beta^{-1}b^{-\beta}\mu(A)\mu(B)
+O(b^{-(2\beta-\beta_+)})\cos tb\,db
\\ & = \Re c_\beta^{-1}\int_0^\omega b^{-\beta}\cos b\,db\,\mu(A)\mu(B)+
O(t^{-(\beta_+-\beta)}\omega^{1-2\beta+\beta_+}).
\end{align*}
Finally, a calculation (see for example~\cite[Proposition~9.5]{MT12}) shows that
$\int_0^\omega b^{-\beta}\cos b\,db=\Gamma(1-\beta)\sin (\beta\pi/2)+O(\omega^{-\beta})$.
Hence the result follows with $\omega=t^\kappa$ for any $\kappa>0$.
\end{proof}

\begin{lemma} \label{lem:I2A}  
$\lim_{t\to\infty}m(t)I_2(t,t^\kappa)=0$ 
for all $\kappa>0$ sufficiently large.
\end{lemma}

\begin{proof}  
We use the same decomposition 
$m(t)|I_2(t,\omega)|  \le a^{-1}J_1(t,\omega)+ \pi a^{-2}J_2(t,\omega)+ a^{-1}J_3(t,\omega)$ as in the proof of Lemma~\ref{lem:I2}.
By Proposition~\ref{prop:psi}, we still have 
$J_1(t,\omega)\ll \omega^{-(\beta-\eps)}$ and
$J_2(t,\omega)\ll t^{-\beta}$.
By Lemma~\ref{lem:TcontA} with $h=\pi/t$,
\[
J_3(t,\omega)\ll  t^{1-\beta}t^{-(\beta-\eps)}\int_{\omega/t}^\infty b^{-2\beta}\,db\ll t^{\eps} \omega^{-(2\beta-1)},
\]
for any choice of $\eps>0$.
Now take $\omega=t^\kappa$ with $\eps<\kappa(2\beta-1)$.
\end{proof}

\section{Proof of the local limit theorem with error term}
\label{sec:LLT}

In this section, we prove Theorem~\ref{thm:LLT}.
The proof combines results from Section~\ref{sec:BMT} with arguments from~\cite{Stone65}.
(A related argument~\cite[Theorem 6.3]{AaronsonDenker01} based on~\cite{Breiman}
gives a similar conclusion but without the error term.)

For ease of exposition, we assume hypotheses~(H) and~(S)(ii) throughout.  However, Lemma~\ref{lem:Tcont} is not required in this section, so we can just as well use hypothesis~(A) instead of hypothesis~(H) by Proposition~\ref{prop:psi}.
Recall that 
$q_\beta(t)=\frac{1}{2\pi}\int_{-\infty}^\infty e^{ibt}e^{-c_\beta|b|^\beta}\,db$ where
$c_\beta=i\int_0^\infty e^{-i\sigma}\sigma^{-\beta}\,d\sigma$.

In Section~\ref{sec:steps}, we made use of the family of kernels $g_a(b)=a^{-1}g(b/a)$ with Fourier transforms $\hat g_a(x)=\hat g(ax)$, where
\[
g(b)=\begin{cases} 1-|b|, &  |b|\leq 1 \\ \hphantom{Y} 0, & |b|>1
\end{cases} \quad \text{and}\quad
\hat g(x)=\frac{2(1-\cos x)}{x^2}.
\]
Since the current section closely follows~\cite{Stone65} which uses slightly different conventions, we now use
$k_a(b)=g(ab)$ with transforms $\hat k_a(x)=\frac{1}{2\pi}a^{-1}\hat g(b/a)$.
(In~\cite{Stone65}, $\hat k_a$ is called $K_a$.)

Let
\[
\mu_n(I)=\mu(y\in A\cap F^{-n}B:\tau_n(y)\in I),
\]
and define
\[
V_n(t,h,a)=\int_{-\infty}^\infty \hat k_a(t-t')\mu_n([d_nt',d_n(t'+h)])\,dt'.
\]

\begin{lemma} \label{lem:Stone}
Let $L>0$.  Then
\[
V_n(t,h,a)=h\big\{q_\beta(t)\mu(A)\mu(B)+e(n,h,a,t)\big\}\quad\text{for $a\ge (Ld_n)^{-1}$},
\]
where $e(n,h,a,t)\to0$ 
as $n\to\infty$, $h\to0$ and $a\to0$, uniformly in $t\in\R$.
\end{lemma}

\begin{proof}
In fact, we show that
\[
|V_n(t,h,a)-hq_\beta(t)\mu(A)\mu(B)|\le {\rm const.}\, h\{e_1(n)+e_2(h)+e_3(a)\}
\]
where $\lim_{n\to\infty}e_1(n)=
\lim_{h\to0}e_2(h)=
\lim_{a\to0}e_3(a)=0$.

As in Section~\ref{sec:BMT}, we write
$\hat R(ib)=\lambda(ib)P(ib)+\tilde Q(b)$ for $|b|\le\delta$,
where $\tilde Q(b)=R(ib)Q(ib)$.
Then
\begin{align}
\label{eq:RPQ}
\hat R(ib)^n=\lambda(ib)^n P(0) +\lambda(ib)^n(P(ib)-P(0)) +\tilde Q(b)^n.
\end{align}
Moreover, there exist constants $C>0$, $\gamma>1-\beta$, $\alpha_1\in(0,1)$,
where
\begin{align} \label{eq:PQ}
\|P(ib)-P(0)\|_{\cB\to L^1}\leq C|b|^\gamma, \quad
\|\tilde Q(b)^n\|_\cB\leq C\alpha_1^n, \quad\text{for all $|b|\le\delta$, $n\ge1$}.
\end{align}

Also, 
we can choose $C>0$, $\alpha_1\in(0,1)$ so that
\begin{align} \label{eq:blarge}
\|\hat R(ib)^n\|_\cB\le C\alpha_1^n\quad\text{for all $b\in[\delta,L]$, $n\ge 1$}.
\end{align}
(Such an estimate for fixed $b>0$ holds by~(S)(ii).  The uniform estimate follows from~\cite[Corollary~2, part 2]{KellerLiverani99}.)

By Corollary~\ref{cor:lambda}(a),
$1-\lambda(ib)\sim c_\beta\ell(1/|b|)b^\beta$.
Hence
\begin{align}
\label{eq:lamasy}
\lambda(ib)\sim e^{-c_\beta \ell(1/|b|)|b|^\beta}\;\text{as $b\to0$},\qquad \lim_{n\to\infty}\lambda(id_n^{-1}b)^n= e^{-c_\beta |b|^\beta}.
\end{align}
Let $\beta'\in(0,\beta)$.
By~\eqref{eq:lamasy} and Potter's bounds, for each fixed $n$,  
there exists $C_1(n),C_2(n)>0$ such that
$|\lambda(id_n^{-1}b)|^n\leq C_1(n) e^{-C_2(n) |b|^{\beta'}}$ for all $|b|\le\delta d_n$.
Also, there exists $n_0\ge1$ such that
$|\lambda(id_n^{-1}b)|^n\leq 2 e^{-c_\beta |b|^{\beta}}$ for all $|b|\le\delta d_n$, $n\ge n_0$.  
Hence
there exists $C_1,C_2>0$ such that
\begin{align}
\label{eq:lamasy2}
|\lambda(id_n^{-1}b)|^n\leq C_1 e^{-C_2 |b|^{\beta'}} \quad\text{for all $|b|\le\delta d_n$, $n\ge1$}.
\end{align}

Now $\hat k_a(t)=\frac{1}{2\pi}\int_\R e^{-ibt}k_a(b)\,db$ and hence
\begin{align*}
V_n(t,h,a) & = \frac{1}{2\pi}\int_{-\infty}^\infty  \int_{-\infty}^\infty e^{-ib(t-t')}k_a(b)\,db\int_{A\cap F^{-n}B}1_{\{\tau_n\in [d_nt',d_n(t'+h)]\}}\,d\mu\,dt'
\\ & = \frac{1}{2\pi}\int_{|b|\le a^{-1}}  e^{-ibt}k_a(b)\int_{A\cap F^{-n}B}\int_{d_n^{-1}\tau_n-h}^{d_n^{-1}\tau_n}e^{ibt'}\,dt'\,d\mu\,db
\\ & =\frac{1}{2\pi}\int_{|b|\le a^{-1}} e^{-itb}k_a(b)\,(1-e^{-ihb})\,(ib)^{-1} \int_{A\cap F^{-n}B} e^{id_n^{-1}b\tau_n}\,d\mu\,db \\
& =\frac{h}{2\pi}\int_{|b|\le a^{-1}} e^{-itb}G(b,h,a) \int_B \hat R(id_n^{-1}b)^n1_A\,d\mu\,db,
\end{align*}
where $G(b,h,a)=k_a(b)\,(1-e^{-ihb})\,(ihb)^{-1}$.

Note that $|G(b,h,a)|\le 1$.  
Using~\eqref{eq:blarge} and that $a\ge (Ld_n)^{-1}$,
\begin{align*}
\Big|\int_{\delta d_n\le |b|\le a^{-1}} e^{-itb}G(b,h,a) \int_B &  \hat R(id_n^{-1}b)^n1_A\,d\mu\,db\Big|
 \le \|1_A\|_{\cB}\int_{\delta d_n\le |b|\le Ld_n}\|\hat R(id_n^{-1}b)^n\|_\cB\,db
\\ & = \|1_A\|_\cB\, d_n \int_{\delta \le |b|\le L}\|\hat R(ib)^n\|_\cB\,db
\le C\|1_A\|_\cB\, d_n\alpha_1^n.
\end{align*}
Hence this term can be incorporated into $e_1(n)$.

It remains to analyse
\[
\frac{h}{2\pi}\int_{|b|\le \delta d_n} e^{-itb}G(b,h,a) \int_B \hat R(id_n^{-1}b)^n1_A\,d\mu\,db
=\frac{h}{2\pi}(I_1+I_2+I_3),
\]
where by~\eqref{eq:RPQ},
\begin{align*}
I_1 & =\int_{|b|\le \delta d_n} e^{-itb}G(b,h,a) \int_B \lambda(id_n^{-1}b)^nP(0)1_A\,d\mu\,db, \\
I_2 & =\int_{|b|\le \delta d_n} e^{-itb}G(b,h,a) \int_B \lambda(id_n^{-1}b)^n(P(id_n^{-1}b)-P(0))1_A\,d\mu\,db, \\
I_3 & =\int_{|b|\le \delta d_n} e^{-itb}G(b,h,a) \int_B \tilde Q(d_n^{-1}b)^n1_A\,d\mu\,db.
\end{align*}
By~\eqref{eq:PQ} and~\eqref{eq:lamasy2},
\[
|I_2|\le \int_{|b|\le\delta d_n} C_1e^{-C_2|b|^{\beta'}}C|d_n^{-1}b|^\gamma\|1_A\|_\cB\,db
\le CC_1\|1_A\|_\cB\,  d_n^{-\gamma}\int_{-\infty}^\infty |b|^\gamma e^{-C_2|b|^{\beta'}}db \ll d_n^{-\gamma},
\]
and
\[
|I_3|\le d_n\int_{|b|\le\delta}C\alpha_1^n\|1_A\|_\cB\,db
\ll d_n\alpha_1^n.
\]
Again, these terms can be incorporated into $e_1(n)$.

This leaves the term $I_1=I_1'\mu(A)\mu(B)$ where
$I_1' =\int_{|b|\le \delta d_n} e^{-itb}G(b,h,a) \lambda(id_n^{-1}b)^n\,db$.
Write $I_1'=J_1+J_2+J_3$ where
\begin{align*}
J_1  & =\int_{|b|\le \delta d_n} e^{-itb} k_a(b)\big\{(1-e^{-ihb})(ihb)^{-1}-1\big\} \lambda(id_n^{-1}b)^n\,db,
\\
J_2  & =\int_{|b|\le \delta d_n} e^{-itb} (k_a(b)-1)\lambda(id_n^{-1}b)^n\,db,
\\
J_3  & =\int_{|b|\le \delta d_n} e^{-itb} \lambda(id_n^{-1}b)^n\,db.
\end{align*}
Since $|(1-e^{-ihb})(ihb)^{-1}-1|\le \frac12 h|b|$ it follows from~\eqref{eq:lamasy2} that
\[
|J_1|\le h\int_{-\infty}^\infty C_1e^{-C_2|b|^{\beta'}}|b|\,db\ll h.
\]
Also,
\[
|J_2|\le \int_{-\infty}^\infty |k_a(b)-1|C_1e^{-C_2|b|^{\beta'}}\,db,
\]
which converges to zero by the dominated convergence theorem as $a\to0$.
These are the sole contributions to $e_2$ and $e_3$ respectively.

Finally, 
\[
|J_3-2\pi q_\beta(t)|  \le 
\int_{|b|\le\delta d_n}|\lambda(id_n^{-1}b)^n- e^{-c_\beta |b|^\beta}|\,db+
\int_{|b|\ge\delta d_n}
e^{-c_\beta |b|^\beta}\,db,
\]
which converges to zero by~\eqref{eq:lamasy},~\eqref{eq:lamasy2} and the dominated convergence theorem as $n\to\infty$.
\end{proof}

\begin{lemma} \label{lem:LLT'}
Let $\eps>0$ and $L>0$. There exists $n_0\ge1$ and $h_0>0$ such that
\[
h(q_\beta(t)\mu(A)\mu(B)-\eps)\le \mu_n([d_nt,d_n(t+h)])\le h(q_\beta(t)\mu(A)\mu(B)+\eps),
\]
for all $n\ge n_0$, $h\in[(Ld_n)^{-1},h_0]$, $t\in\R$.
\end{lemma}

\begin{proof}
Let $\tilde q_\beta=q_\beta\mu(A)\mu(B)$.
Since $q_\beta$ is the Fourier transform of an $L^1$ function, 
$\tilde q_\beta$ is uniformly continuous and bounded.
Let $q_\infty=|\tilde q_\beta|_\infty$ and 
choose $h_1\in(0,1)$ such that
$|\tilde q_\beta(t)-\tilde q_\beta(t')|\le \frac14\eps$ whenever $|t-t'|\le h_1$.

For $\eps_1>0$, set 
$\eps_2=\int_{|x|>1/\eps_1}\hat k_1(x)\,dx$.  We choose $\eps_1\in(0,\frac16)$ sufficiently small that
\begin{align} \label{eq:q}
\SMALL
(q_\infty+2\eps_1 q_\infty+\frac12\eps)(1-\eps_2)^{-1}-q_\infty\le\eps,
\qquad
2\eps_1 q_\infty +\eps_2(q_\infty+\eps)\le\frac12 \eps.
\end{align} 

By Lemma~\ref{lem:Stone},
there exists $n_0\ge1$ and $h_0\in(0,h_1)$ such that
for all $n\ge n_0$, $h\in[(Ld_n)^{-1},h_0]$, $t\in\R$,
\begin{align} \label{eq:V1}
\SMALL
V_n(t-\eps_1 h,h(1+2\eps_1),\eps_1^2h)
& \SMALL \le
h(1+2\eps_1)\tilde q_\beta(t-\eps_1 h)+\frac16 \eps h
\\ & \SMALL \le
h(1+2\eps_1)(\tilde q_\beta(t)+\frac14\eps)+\frac16 \eps h
\le
h(\tilde q_\beta(t)+2\eps_1 q_\infty+\frac12\eps), \nonumber
\end{align}
where we used the constraint $\eps_1\le\frac16$.
Also, we can ensure that
\begin{align} \label{eq:V2}
\SMALL
V_n(t+\eps_1 h,h(1-2\eps_1 ),\eps_1^2h)
& \SMALL \ge
h(1-2\eps_1)\tilde q_\beta(t+\eps_1 h)-\frac14 \eps h
\\ & \SMALL \ge
h(1-2\eps_1)(\tilde q_\beta(t)-\frac14\eps)-\frac14 \eps h
\ge
h(\tilde q_\beta(t)-2\eps_1 q_\infty-\frac12\eps). \nonumber
\end{align}

Now, for $|t'|\le \eps_1 h$,
\begin{align*}
\mu_n([d_n(t+\eps_1 h-t')  \SMALL,d_n(t-\eps_1 h-t'+h])
 &  \le \mu_n([d_nt,d_n(t+h)])
\\ & \le \mu_n([d_n(t-\eps_1 h-t'),d_n(t+\eps_1 h-t'+h]).
\end{align*}
Also $\int_{-\infty}^\infty \hat k_1\,dx=1$, so
\[
1-\eps_2=\int_{|x|\le 1/\eps_1}\hat k_1(x)\,dx=
\eps_1^2h\int_{|x|\le 1/\eps_1}\hat k_{\eps_1^2h}(\eps_1^2hx)\,dx=
\int_{|x|\le \eps_1 h}\hat k_{\eps_1^2h}(x)\,dx.
\]
Hence
\begin{align*}
& \SMALL V_n(t-\eps_1 h,h(1+2\eps_1),\eps_1^2h)
 \\
 & \qquad = \SMALL
\int_{-\infty}^\infty \hat k_{\eps_1^2h}(t')\mu_n([d_n(t-\eps_1 h-t'),d_n(t+\eps_1 h-t'+h)])\,dt'
\\
 & \qquad \ge \SMALL
\int_{|t'|\le \eps_1 h} \hat k_{\eps_1^2h}(t')\mu_n([d_n(t-\eps_1 h-t'),d_n(t+\eps_1 h-t'+h)])\,dt'
 \\ & \qquad \ge \SMALL
\int_{|t'|\le \eps_1 h} \hat k_{\eps_1^2h}(t')\mu_n([d_nt,d_n(t+h)])\,dt'
=(1-\eps_2) \mu_n([d_nt,d_n(t+h)]).
\end{align*}
By~\eqref{eq:q} and~\eqref{eq:V1},
\begin{align*}
\mu_n([d_nt,d_n(t+h)]) & \SMALL \le(1-\eps_2)^{-1}V_n(t-\eps_1 h,h(1+2\eps_1),\eps_1^2h)
\\ & \SMALL \le h(\tilde q_\beta(t)+2\eps_1 q_\infty+\frac12\eps)(1-\eps_2)^{-1}
\le h(\tilde q_\beta(t)+\eps).
\end{align*}

Arguing similarly, and exploiting the last estimate for $\mu_n([d_nt,d_n(t+h)])$,
\begin{align*}
& \SMALL V_n(t+\eps_1 h,h(1-2\eps_1),\eps_1^2h)
 \\
 & \qquad \le \SMALL
\int_{|t'|\le \eps_1 h} \hat k_{\eps_1^2h}(t')\mu_n([d_n(t+\eps_1 h-t'),d_n(t-\eps_1 h-t'+h)])\,dt' 
\\ & \qquad\qquad\qquad\qquad + \int_{|t'|\ge \eps_1 h} \hat k_{\eps_1^2h}(t')h(q_\infty+\eps)\,dt'
 \\ & \qquad \le \SMALL
\mu_n([d_nt,d_n(t+h)]) + \eps_2 h (q_\infty+\eps).
\end{align*}
By~\eqref{eq:q} and~\eqref{eq:V2},
\begin{align*}
\mu_n([d_nt,d_n(t+h)]) & \SMALL \ge V_n(t+\eps_1 h,h(1-2\eps_1),\eps_1^2h) - \eps_2 h(q_\infty+\eps)
\\ & \SMALL \ge h((\tilde q_\beta(t)-2\eps_1 q_\infty-\frac12\eps - \eps_2(q_\infty+\eps)) 
\ge h(\tilde q_\beta(t)-\eps).
\end{align*}
This completes the proof.
\end{proof}

\begin{pfof}{Theorem~\ref{thm:LLT}}
After a change of variables, Lemma~\ref{lem:LLT'} reads as follows:

Let $\eps>0$ and $L>0$. There exists $n_0\ge1$ and $h_0>0$ such that
\begin{align}\label{eq:LLT}
\sup_{t\in\R}d_n\Big|\mu_n([t,t+h])-\frac{h}{d_n}q_\beta(d_n^{-1}t)\mu(A)\mu(B)\Big|\le h\eps,
\end{align}
for all $n\ge n_0$, $h\in[L^{-1},d_nh_0]$.

Fix $h>0$ and define $e_n=\sup_{t\in\R}d_n\big|\mu_n([t,t+h])-\frac{h}{d_n}q_\beta(d_n^{-1}t)\mu(A)\mu(B)\big|$.
We must show that $\lim_{n\to\infty}e_n=0$.

Let $L= 1/h$.  
By~\eqref{eq:LLT}, for any $\eps>0$ there exists $n_0\ge1$, $h_0>0$, 
such that $e_n\le h\eps$ for all $n\ge n_0$ subject to the constraint
$d_nh_0\ge h$.   Since $d_n\to\infty$, there exists $n_1\ge n_0$ such that 
$d_nh_0\ge h$ for all $n\ge n_1$. Hence 
$e_n\le h\eps$ for all $n\ge n_1$ as required.~
\end{pfof}

\section{Proof of Theorems~\ref{thm:inf} and~\ref{thm:average}}
\label{sec:pfthmaver}

In this section, we prove Theorem~\ref{thm:inf} by establishing separately an upper bound
(Corollary~\ref{cor:infupper}) and a lower bound (Corollary~\ref{cor:inflower}).
In the process of obtaining the upper bound, we prove Theorem~\ref{thm:average}.

For ease of exposition, we assume hypothesis~(H) throughout.  Again, Lemma~\ref{lem:Tcont} is not required in this section, so we can just as well use hypothesis~(A) by Proposition~\ref{prop:psi}.

\subsection{Upper bound for $\liminf$}
\label{sec:infupper}

In this subsection, we only require hypothesis~(H)
with $s\in\R^+$ in (H)(ii).
A simplified version of the argument used in the proof of Lemma~\ref{lem:T}
can be used to obtain

\begin{prop} \label{prop:real}
Assume the setting of Theorem~\ref{thm:average} with $\beta\in[0,1]$. 
For $\sigma>0$,
\[
\hat T(\sigma)= 
{D_\beta}'\,\tilde\ell(1/\sigma)^{-1}\sigma^{-\beta}(P(0)+E(\sigma)), 
\]
where ${D_\beta}'= \Gamma(1-\beta)^{-1}$ for $\beta\in(0,1)$ and
${D_0}'={D_1}'=1$, and
$E(\sigma)$ is a family of operators satisfying $\lim_{\sigma\to 0}\|E(\sigma)\|_{\cB\to L^1}=0$. \qed
\end{prop}

We can now complete

\begin{pfof}{Theorem~\ref{thm:average}}
For $n\ge0$, the real Laplace transform of the distribution $G_n(x)=\mu(\tau_n(y)\leq x,\,y\in A\cap F^{-n}B)$
is given by $\int_Y 1_A\,1_B\circ F^n\,e^{-\sigma\tau_n}\, d\mu=\int_B \hat R(e^{-\sigma})^n 1_A\,d\mu$.
 Hence,
\[
 \int_{-\infty}^{\infty} e^{-\sigma t}\, d U_{A,B}(t)
=\sum_{n=0}^\infty \int_B \hat R(e^{-\sigma})^n 1_A\,d\mu
=\int_B \hat T(e^{-\sigma}) 1_A\,d\mu.
\]
The conclusion follows  from
Proposition~\ref{prop:real} by the continuous time version of
Karamata's Tauberian Theorem~\cite[Theorem 1.7.1]{BGT}.
\end{pfof}

\begin{lemma} \label{lem:9}
Assume the setting of Theorem~\ref{thm:average} with $\beta\in(0,1]$. 
Let $z:[0,\infty)\to[0,\infty)$ be integrable.   Then
\[
\liminf_{t\to\infty} m(t)\int_0^t z(t-y)\,dU_{A,B}(y)\le d_\beta \mu(A)\mu(B) \int_0^\infty z\,dx.
\]
\end{lemma}

\begin{proof}  This is proved in the same way as~\cite[Lemma~9]{Erickson70}
using Theorem~\ref{thm:average}.
\end{proof}

\begin{cor}
\label{cor:infupper}
Assume the setting of Theorem~\ref{thm:average} with $\beta\in(0,1]$. 
Then for any $h>0$, 
\[
\liminf_{t\to\infty}m(t)(U_{A,B}(t+h)-U_{A,B}(t))\leq d_{\beta}\mu(A)\mu(B)h.
\]
\end{cor}

\begin{proof}
Let $z=1_{[0,h]}$.  By Lemma~\ref{lem:9},
\begin{align*}
\liminf_{t\to\infty} m(t)(U_{A,B}(t+h)-U_{A,B}(t)) & =\liminf_{t\to\infty} m(t+h)\int_0^{t+h} z(t+h-y)\,dU_{A.B}(y) \\
& \le d_\beta \mu(A)\mu(B)\int_0^\infty z\,dx
= d_\beta \mu(A)\mu(B)h,
\end{align*}
as required.
\end{proof}

\subsection{Lower bound for $\liminf$}
\label{sec:inflower}


\begin{cor}
\label{cor:inflower}Assume the setting of Theorem~\ref{thm:inf}. Then for any $h>0$, 
\[
\liminf_{t\to\infty}m(t)(U_{A,B}(t+h)-U_{A,B}(t))\geq d_{\beta}\mu(A)\mu(B)h.
\]
\end{cor}

\begin{proof}
Let $m\ge k\ge0$.
By~\eqref{eq:U} and Theorem~\ref{thm:LLT}, 
\begin{align*}
U_{A,B}(t+h)-U_{A,B}(t)&\geq \sum_{n=k}^m \mu(y\in A\cap F^{-n}B:\tau_n(y)\in[t,t+h] )\\
&=\sum_{n=k}^m \frac{h}{d_n}q_\beta(t/d_n)\mu(A)\mu(B)+E_{k,m},
\end{align*}
where $E_{k,m}=\sum_{n=k}^m e_n/d_n$.

Let $\kappa\in(1,1/\beta)$.  Then
$d_n^{-1}=O(n^{-\kappa})$ and 
$E_{k,m}=O(\sup_{n\ge k}|e_n|)\to0$ as $k\to\infty$.

Choosing $k=[C_1t^\beta/\ell(t)]$ and $m=[C_2t^\beta/\ell(t)]$, for fixed $C_2> C_1>0$ and arguing word for word as
in~\cite[Proof of eq.~(7.2)]{Erickson70}, we obtain
\[
\liminf_{t\to\infty}m(t)(U_{A,B}(t+h)-U_{A,B}(t))\geq\mu(A)\mu(B)\int_{C_1}^{C_2}x^{-1/\beta}q_\beta(x^{-1/\beta})\, dx.
\]
Now let $C_1\to 0$ and $C_2\to \infty$ and use that
$\int_{0}^{\infty}x^{-1/\beta}q_\beta(x^{-1/\beta})\, dx=d_\beta$.
\end{proof}

\section{General class of observables}
\label{sec:obs}

In this section, we extend mixing for semiflows, Corollary~\ref{cor:mixing}, to cover more general classes of observables.  As well as being of interest in its own right, this is useful for the extension to flows in Section~\ref{sec:flow}.

Throughout, we suppose that we are in the setting of Corollary~\ref{cor:mixing}; in particular $\beta\in(\frac12,1]$ and hypotheses (H) and (S)(i) hold.
We also suppose from now on that $Y$ is a metric space with inner regular\footnote{$\mu$ is inner regular if $\mu(A)=\sup\mu(K: K\subset A,\;\text{$A$ compact}\}$ for all open sets $A\subset Y$.}
Borel probability measure $\mu$ and that $F$ and $\tau$ are almost everywhere continuous.
It is well-known that mixing for infinite measure system is not a measure-theoretic property~\cite{HajianKakutani64,Krickeberg67} and that care needs to be taken with the class of observables.  Here we follow Krickeberg~\cite{Krickeberg67}.
As a special case of the general theory, we prove the following result:

\begin{thm} \label{thm:H}
Define $H_n=\{(y,u)\in Y\times[0,\infty):\tau(y)-n\le u\le \tau(y)\}$, $n\ge1$.
Then
\begin{align} \label{eq:H}
\lim_{t\to\infty} m(t)\int_{Y^\tau} v\,w\circ F_t\,d\mu^\tau  =d_\beta
\int_{Y^\tau} v\,d\mu^\tau
\int_{Y^\tau} w\,d\mu^\tau
\end{align} 
for all bounded and almost everywhere continuous functions $v:Y^\tau\to\R$ supported in $H_n$ for some $n$, and all $w\in L^1(Y^\tau)$.
\end{thm}

Note that this includes all bounded almost everywhere continuous observables $v$ supported in a set of the form $A\times[a_1,a_2]\subset Y^\tau$
where $A\subset Y$, $0<a_1<a_2\le\inf_A\tau$ and $\sup_A\tau<\infty$.  For the results on flows in Section~\ref{sec:flow} we require the more general class of observables in Theorem~\ref{thm:H}.

In the remainder of this section, we prove a more general result along the lines of~\cite{Krickeberg67} and use this to prove Theorem~\ref{thm:H}.

Let $\cC$ be a collection of measurable subsets $A\subset Y$ with
$1_A\in\cB$ such that
\begin{itemize}
\parskip=-2pt
\item[(i)] $\mu(\partial A)=0$ for all $A\in\cC$,
\item[(ii)] $A_1\cap A_2\in\cC$ for all $A_1,A_2\in\cC$,
\item[(iii)] $\cC$ is a basis for the topology on $Y$.
\end{itemize}
In practice, we can often take $\cC$ to consist of 
all measurable sets $A\subset Y$ with $1_A\in\cB$ and $\mu(\partial A)=0$.
This is the case for the examples in Section~\ref{sec:ex}.

\begin{prop} \label{prop:obs}
Let $\cC'=\{A\times[a_1,a_2]\subset Y^\tau: A\in\cC\}$.
Let $\cD$ be the ring generated by $\cC'$
and let $H\in\cD$.
Then~\eqref{eq:H} holds
for all bounded and almost everywhere continuous functions \mbox{$v:Y^\tau\to\R$} supported in $H$, and all $w\in L^1(Y^\tau)$.
\end{prop}

\begin{proof}
It is immediate that conditions~(i)--(iii) for $\cC$ are inherited by 
the collection $\cC'$ of subsets of $Y^\tau$ 
(with $\mu$ replaced by $\mu^\tau$).  

Write $q(t)=d_\beta^{-1}m(t)$.  By Corollary~\ref{cor:mixing}, 
\begin{align} \label{eq:mix}
\lim_{t\to\infty} q(t)\mu^\tau(A\cap F_t^{-1}B)  =\mu^\tau(A)\mu^\tau(B),
\end{align} 
for all $A\in\cC'$ and all measurable rectangles $B\subset Y^\tau$.
The argument now proceeds as in~\cite[Section 2]{Krickeberg67}.
We provide the details for completeness.

\noindent{\bf Step 1:}
Let $B\subset Y^\tau$ be a measurable rectangle.  
Then~\eqref{eq:mix} holds for all $A\in\cC'$ and hence (using condition~(ii)) for all finite unions and differences of elements of $\cC'$.  In other words,~\eqref{eq:mix} holds for all $A\in\cD$.

\noindent{\bf Step 2:}
Let $B\subset Y^\tau$ be a measurable rectangle.
Recall that $H\in\cD$ and
let $A\subset H$ such that $\mu^\tau(\partial A)=0$.  
Suppose that $K\subset \Int A$ is compact.
Since $\cC'$ is a basis and $\cD$ is stable under finite unions, 
there exists $D\in\cD$ such that $K\subset D\subset\Int A$.  
Using also the inner regularity of $\mu^\tau$,
\begin{align*}
\mu^\tau(A)=\mu^\tau(\Int A)
& =\sup\{\mu^\tau(K):K\subset\Int A,\;K\,\text{compact}\}
\\ & = \sup\{\mu^\tau(D):D\subset A,\;D\in\cD\}.
\end{align*}
Similarly, 
\[
\mu^\tau(H\setminus A)=\sup\{\mu^\tau(D):D\subset H\setminus A,\;D\in\cD\}
=\sup\{\mu^\tau(H\setminus D):D\supset A,\;D\in\cD\},
\]
so
$\mu^\tau(A)=\inf\{\mu^\tau(D):D\supset A,\;D\in\cD\}$.
Hence for any $\eps>0$, there exist $D_1,D_2\in\cD$ such that 
$D_1\subset A\subset D_2$ and $\mu^\tau(D_2)-\mu^\tau(D_1)<\eps$.
Since~\eqref{eq:mix} holds for $D_1$ and $D_2$ and
  $\mu^\tau(D_1\cap F_t^{-1}B)\le
  \mu^\tau(A\cap F_t^{-1}B)
  \le \mu^\tau(D_2\cap F_t^{-1}B)$,
\begin{align*}
(\mu^\tau(A)-\eps)\mu^\tau(B) & \le 
\mu^\tau(D_1)\mu^\tau(B)
  \le\liminf_{t\to\infty}q(t)\mu^\tau(A\cap F_t^{-1}B)
\\ & \le\limsup_{t\to\infty}q(t)\mu^\tau(A\cap F_t^{-1}B)
\le \mu^\tau(D_2)\mu^\tau(B)
\le (\mu^\tau(A)+\eps)\mu^\tau(B).
\end{align*}
As $\eps$ is arbitrary, we have verified that~\eqref{eq:mix} holds for all
$A\subset H$ with $\mu^\tau(\partial A)=0$.
In other words,
$\lim_{t\to\infty}q(t)\int_{Y^\tau}v\, 1_B\circ F_t\,d\mu^\tau=\int_{Y^\tau}v\,d\mu^\tau\,\mu^\tau(B)$ where $v=1_A$.
This extends to all finite linear combinations $v=\sum c_j1_{A_j}$
by linearity.  We will refer to such functions $v$ as step functions.

\noindent{\bf Step 3:}
Let $B\subset Y^H$ be a measurable rectangle and
suppose that $v$ is as in the statement of the proposition.
We claim that for any $\eps>0$ there exist step functions $v_1$ and $v_2$ such that
$v_1\le v\le v_2$ and $\int_{Y^\tau}v_2\,d\mu^\tau -  \int_{Y^\tau}v_1\,d\mu^\tau<\eps$.
Then
\begin{align*}
& \Big(\int_{Y^\tau}v\,d\mu^\tau-\eps\Big)\mu^\tau(B)
 \le  \int_{Y^\tau}v_1\,d\mu^\tau\,\mu^\tau(B)
\le \liminf_{t\to\infty}q(t)\int_{Y^\tau}v\, 1_B\circ F_t\,d\mu^\tau
\\ &\qquad \le \limsup_{t\to\infty}q(t)\int_{Y^\tau}v\, 1_B\circ F_t\,d\mu^\tau
\le \int_{Y^\tau}v_2\,d\mu^\tau\,\mu^\tau(B)
\le \Big(\int_{Y^\tau}v\,d\mu^\tau+\eps\Big)\mu^\tau(B).
\end{align*}
Hence~\eqref{eq:H} holds for all $v$ of the desired form and all
indicator functions $w=1_B$ where $B$ is a measurable rectangle.

To prove the claim, let $\delta>0$ such that $\delta(\mu^\tau(Y)+2|v|_\infty)<\eps/2$ and let $I$ be a closed interval covering the image of $v$.
We can write $I$ as a finite union of closed intervals $I_1,\dots,I_N$ with $\diam I_j<\delta$ intersecting only at endpoints.

Let $A_j=v^{-1}(I_j)$ and define $Z$ to be the set of discontinuity points
of $v$.  Then
$\partial A_j\subset Z\cup v^{-1}(\partial I_j)$ for all $j$.
Hence $\mu^\tau(\partial A_j)\le \mu^\tau(v^{-1}(\partial I_j))$.

Also, there are at most countably many $x_k\in\R$ such that $\mu^\tau(v^{-1}(x_k))>0$.  We can modify the intervals $I_j$ slightly so that $x_k\not\in\partial I_j$ for all $j,k$.  This ensures that
$\mu^\tau(\partial A_j)=0$ for all $j$.

As in Step~2, it follows from  inner regularity of $\mu^\tau$ that
for each $j$ there exists $D_j\in\cD$ with $D_j\subset A_j$ such that
$\mu^\tau(A_j\setminus D_j)<\delta/N$.
Now define
\[
\SMALL v_1=\sum \inf_{D_j}v\,1_{D_j} \;+\;\inf_Y v\,1_{H\setminus \bigcup D_j},
\qquad
v_2=\sum \sup_{D_j}v\,1_{D_j} \;+\;\sup_Y v\,1_{H\setminus \bigcup D_j}.
\]
Then $v_1\le v\le v_2$.  Also,
\begin{align*}
\int_{Y^\tau} v_2\,d\mu^\tau-\int_{Y^\tau} v\,d\mu^\tau & \le \SMALL \sum\mu^\tau(D_j)(\sup_{D_j}v-\inf_{D_j}v)\,+\,
2\mu^\tau(H\setminus\bigcup D_j)|v|_\infty
\\ & \le \mu^\tau(Y)\delta + 2|v|_\infty \sum \mu^\tau(A_j\setminus D_j)
< \delta(\mu^\tau(Y)+2|v|_\infty)<\eps/2.
\end{align*}
Similarly, $\int_{Y^\tau} v\,d\mu^\tau-\int_{Y^\tau} v_1\,d\mu^\tau<\eps/2$ verifying  the claim.

\noindent{\bf Step 4:}
To prove the general result, suppose without loss that $v\ge0$ and let $w\in L^1(Y^\tau)$.
By a more standard approximation argument than the one in Step~3, there exist simple functions
$w_1$ and $w_2$ such that
$w_1\le w\le w_2$ and $\int_{Y^\tau}w_2\,d\mu^\tau -  \int_{Y^\tau}w_1\,d\mu^\tau{<\eps}$.
The result follows.
\end{proof}

\begin{pfof}{Theorem~\ref{thm:H}}
Let $\cC''=\cC'\cup\{E_n,\,n\ge1\}$ where $\cC'$ is the collection of rectangles in Proposition~\ref{prop:obs}
and $E_n=\bigcup_{j=1}^nF_j^{-1}(Y\times[0,1])$.  Let $\cI=\{C\cap E_n:C\in \cC',\,n\ge1\}$
and define $\cC'''=\cC''\cup \cI$.
Then $\cC'''$ is closed under finite intersections, and 
hence conditions~(i)--(iii) are satisfied by the collection $\cC'''$.
We claim that property~\eqref{eq:mix} holds for all $A\in\cC'''$.
Certainly, the sets $E_n$ lie in the ring generated by $\cC'''$, and $H_n\subset E_n$, so the conclusion follows from the approximate argument used to prove Proposition~\ref{prop:obs}.

It remains to verify the claim.  By Corollary~\ref{cor:mixing}, property~\eqref{eq:mix} holds for
all $A\in \cC'$.  By Remark~\ref{rmk:mixing}, this holds also for the sets
$E_n$.
Finally, if $I\in\cI$, then $I$ is contained in one of the rectangles
in $\cC'$ and $\mu^\tau(\partial I)=0$.
Hence $1_I$ is a bounded and almost everywhere continuous function supported in
a rectangle in $\cC'$.
The claim follows from Proposition~\ref{prop:obs}.
\end{pfof}

\section{Mixing for infinite measure flows}
\label{sec:flow}

In this section, we show how mixing for semiflows extends to mixing for flows.  

\subsection{Assumptions and disintegration}
We suppose throughout that $F_t:Y^\tau\to Y^\tau$ is a suspension semiflow
over a map $F:Y\to Y$ with nonintegrable almost everywhere continuous roof function $\tau:Y\to\R^+$ satisfying $\essinf\tau>1$ and $\mu(\tau>t)=\ell(t)t^{-\beta}$, $\beta\in(\frac12,1]$, and we assume that hypotheses (H) and (S)(i) hold.

Let $X=Y\times N$ where $Y$ and $N$ are bounded metric space.
Let $f(y,z)=(Fy,G(y,z))$ where $F:Y\to Y$ and
$G:Y\times N\to N$ are continuous almost everywhere.
The projection $\pi:X\to Y$,  $\pi(y,z)=y$,  defines a semiconjugacy between $f$ and $F$.  
There exists a unique $f$-invariant ergodic probability measure $\mu_X$ on $X$
such that $\pi_*\mu_X=\mu$,
see for instance~\cite[Section~6]{APPV09}.

Define $\tau:X\to\R^+$ by setting $\tau(y,z)=\tau(y)$ and
define the suspension $X^\tau=\{(x,u)\in X\times\R:0\le u\le \tau(x)\}/\sim$
where $(x,\tau(x))\sim(fx,0)$.  The suspension flow
$f_t:X^\tau\to X^\tau$ is given by $f_t(x,u)=(x,u+t)$ computed modulo identifications,
 with ergodic invariant measure
$\mu_X^\tau=\mu_X\times{\rm Lebesgue}$.

Under two additional assumptions~(F1) and (F2) below, we show in Theorem~\ref{thm:flow} that
Corollary~\ref{cor:mixing} for the semiflow $F_t$ applies equally to
the flow $f_t$.

First, we assume contractivity along $N$:
\begin{itemize}
\item[\textbf{(F1)}]
$\lim_{n\to\infty}d(f^n(y,z),f^n(y,z'))=0$ for all $z,z'\in N$ uniformly in $y\in Y$.
\end{itemize}

Recall that $R$ denotes the transfer operator for $F:Y\to Y$.

\begin{prop} \label{prop:disint}
Fix $z_0\in N$.
Suppose $v\in C^0(X)$.  Then the limit
 \[
  \eta_y(v)=\lim_{n\to\infty} (R^nv_n)(y),\qquad v_n(y)=v\circ f^n(y,z_0),
       \]
exists for almost every $y\in Y$ and defines a probability measure supported on $\pi^{-1}(y)$.
Moreover $y\mapsto \eta_y(v)=\int_{\pi^{-1}(y)}v\,d\eta_y$ is integrable and
$\int_X v\,d\mu_X = \int_Y \int_{\pi^{-1}(y)}v\,d\eta_y \ d\mu(y)$.
\end{prop}

\begin{proof}  See for instance~\cite[Proposition~3]{ButterleyMapp}.
\end{proof}

\begin{rmk} \label{rmk:F2}
The proof of~\cite[Proposition~3]{ButterleyMapp} shows that the sequence $R^nv_n$ is Cauchy in $L^\infty(Y)$.  If the metric on $Y$ can be chosen so that $R^nv_n$ is continuous for each $n$, then $\bar v\in C^0(Y)$.  
(In fact, it can often be shown that $\bar v$ is H\"older when $v$ is H\"older~\cite{ButterleyMapp}.)
\end{rmk}

Note that $X^\tau=Y^\tau\times N$.
Given $v\in C^0(X^\tau)$, define
\[
\bar v:Y^\tau\to\R, \qquad
\bar v(y,u)=\int_{x\in\pi^{-1}(y)} v(x,u)\,d\eta_y(x).
\]
Then
\[
\int_{X^\tau} v\,d\mu_X^\tau = \int_{Y^\tau} \bar v(y,u) \, d\mu^\tau(y,u).
\]

We require the additional assumption:
\begin{itemize}
\item[\textbf{(F2)}]
The function $\bar v:Y^\tau\to\R$ is almost everywhere continuous.
\end{itemize}

\begin{rmk} \label{rmk:F22}
If $v$ is uniformly continuous, then 
for any $\eps>0$ there exists $\delta<0$ such that
$|\bar v(y,u)-\bar v(y,u')|<\eps$ for all $(y,u),(y,u')\in Y^\tau$ with $|u-u'|<\delta$.
This combined with Remark~\ref{rmk:F2} shows that condition~(F2) is easily satisfied in practice for a large class of observables $v\in C^0(X^\tau)$.
\end{rmk}

\begin{rmk} \label{rmk:PH}
The set up in this section (skew product $X=Y\times N$, roof function $\tau$ constant in the $N$ direction) is not very restrictive.  Suppose that $T_t:M\to M$ is a smooth flow defined on a Riemannian manifold $M$ and that $\Lambda$ is a partially hyperbolic attractor, so there exists a continuous $DT_t$-invariant splitting $T_\Lambda M=E^s\oplus E^{cu}$ where $E^s$ is uniformly contracting and dominates $E^{cu}$.  By~\cite[Proposition 3.2, Theorem~4.2]{AraujoMapp}, the stable bundle $E^s$ extends to a neighbourhood $U$ of $\Lambda$ and integrates to a $T_t$-invariant collection $\cW^s$ of stable leaves that topologically foliate $U$.

This means that we can choose a topological submanifold $X\subset M$ that is a cross-section to the flow $T_t$ formed as a union of stable leaves, and automatically the roof function $\tau$ is constant along stable leaves.  (This construction has been widely used recently~\cite{ABV16,AraujoM16,AMV15,AGY06}.)  Assuming for convenience the existence of a global chart for $\cW^s$, we obtain a Poincar\'e map $f:X\to X$ where $X=Y\times N$ with $N$ playing the role of the stable direction.
Moreover,
$f$ has the desired skew product form $f(y,z)=(Fy,G(y,z))$,
where $F:Y\to Y$ is defined by
quotienting along the stable leaves, and condition (F1) is automatically satisfied.  
Also (F2) holds by Remark~\ref{rmk:F2}.
Hence our set up holds in its entirety provided $F:Y\to Y$ and $\tau:Y\to\Z^+$ satisfy the required properties.
\end{rmk}

\subsection{The mixing result}

Choose a subset $H$ of $Y^\tau$ as in Proposition~\ref{prop:obs}.

\begin{thm} \label{thm:flow}
Suppose that $\mu(\tau>n)=\ell(n)n^{-\beta}$ where $\beta\in(\frac12,1]$.
Let $v\in C^0(X^\tau)$ be supported in $C\times N$ where $C$ is a closed subset of
$\Int H$.  
Let $w\in C^0(X^\tau)$ be uniformly continuous and supported on a set of finite measure.
Assume that (H), (S1), (F1) and (F2) hold.
Then
 \[
\lim_{t\to\infty} m(t)\int_{X^\tau} v\,w\circ f_t\,d\mu_X^\tau= d_\beta \int_{X^\tau} v\,d\mu_X^\tau\, \int_{X^\tau} w \,d\mu_X^\tau.
\]
\end{thm}

\begin{proof}
Following~\cite{AGY06}, we define $w_s:Y^\tau\to\R$, $s>0$, by setting
\[
w_s(y,u)=\overline{w\circ f_s}=\int_{x\in\pi^{-1}(y)}w\circ f_s(x,u)\,d\eta_y(x).
\]
Note that $\int_{Y^\tau}|w_s|\,d\mu^\tau\le \int_{X^\tau}|w|\circ f_s\,d\mu_X^\tau
=\int_{X^\tau}|w|\,d\mu_X^\tau$  so $w_s\in L^1(Y^\tau)$ for all $s$.

The semiconjugacy $\pi:X\to Y$ extends to a measure-preserving semiconjugacy $\pi^\tau:X^\tau\to Y^\tau$, $\pi^\tau(x,u)=(\pi x,u)$.
Write $m(t)\int_{X^\tau}v\,w\circ f_t\,d\mu_X^\tau=I_1(s,t)+I_2(s,t)$ where
\begin{align*}
I_1(s,t) & = m(t)\int_{X^\tau}v\,w_s\circ\pi^\tau\circ f_{t-s}\,d\mu_X^\tau, \\
 I_2(s,t) & = m(t)\int_{X^\tau}v\,(w\circ f_s-w_s\circ\pi^\tau)\circ f_{t-s}\,d\mu_X^\tau.
\end{align*}

For $t>s$,
\[
I_1(s,t)    =
m(t)\int_{X^\tau}v\,w_s\circ F_{t-s}\circ \pi^\tau\,d\mu_X^\tau
 = m(t)\int_{Y^\tau} \bar v\, w_s\circ F_{t-s}\,d\mu^\tau.
\]
Since $\bar v$ is bounded and almost everywhere continuous, supported in $H$, and
$w_s\in L^1(Y^\tau)$, it follows from Proposition~\ref{prop:obs} that
for all $s>0$,
\[
\lim_{t\to\infty}I_1(s,t)
=d_\beta\int_{Y^\tau}\bar v\,d\mu^\tau\int_{Y^\tau}w_s\,d\mu^\tau
=d_\beta\int_{X^\tau}v\,d\mu_X^\tau\int_{X^\tau}w\,d\mu_X^\tau.
\]

Choose $\psi:Y^\tau\to[0,1]$ continuous such that
$\supp v\subset \supp\psi\times N\subset H\times N$.
Define 
\[
D_s:Y^\tau\to \R, \qquad D_s(y,u)=\diam w\circ f_s((\pi^\tau)^{-1}(y,u)).
\]
Note that $|D_s|\le 2|w|_\infty$ and $\mu^\tau(\supp D_s)\le \mu_X^\tau(f_s^{-1}\supp w)=\mu_X^\tau(\supp w)<\infty$, so $D_s\in L^1(Y^\tau)$.
Also,
$|w\circ f_s(x,u)- w_s\circ \pi^\tau(x,u)| \le D_s\circ \pi^\tau(x,u)$.
Hence for $t>s$,
\[
|I_2(s,t)|  
 \le |v|_\infty\, m(t)\int_{X^\tau} \psi\circ \pi^\tau\, D_s\circ\pi^\tau\circ f_{t-s}\,d\mu_X^\tau   
= |v|_\infty\, m(t)\int_{Y^\tau} \psi \,D_s\circ F_{t-s}\,d\mu_Y^\tau.
\]
Since $\psi\in C^0(Y^\tau)$ is supported in $H$ and $D_s\in L^1(Y^\tau)$, it
again follows from Proposition~\ref{prop:obs} that for all $s>0$,
\[
\limsup_{t\to\infty}I_2(s,t)\le |v|_\infty\, d_\beta \int_{Y^\tau}\psi\,d\mu^\tau\int_{Y^\tau}D_s\,d\mu^\tau.
\]

By uniform continuity of $w$ and (F1), 
$\lim_{s\to\infty}|D_s|_\infty=0$.
Hence $|D_s|_1 \le |D_s|_\infty\,\mu^\tau(\supp D_s)\le |D_s|_\infty\, \mu_X^\tau(\supp w)\to0$ as $s\to\infty$.
This combined with the estimates for $I_1$ and $I_2$ yields the desired result.
\end{proof}

\section{Examples}
\label{sec:ex}

In this section, we demonstrate how the methods in this paper apply to 
the examples described in the introduction.

\subsection{NonMarkovian intermittent semiflows and flows.}
\label{sec:AFU}

Let $f_t:[0,1]^{\tau_0}\to[0,1]^{\tau_0}$ be an intermittent semiflow as in Example~\ref{ex:neutral}.
The first step is to pass from the original suspension semiflow on
$[0,1]^{\tau_0}$ to a suspension of the form
$Y^\tau$ where $(Y,\mu)$ is a probability space and $\tau$ is an nonintegrable roof function.

We take $Y\subset[0,1]$ to be the interval of domain of the rightmost branch of the AFN map $f:[0,1]\to[0,1]$.
Define the first return map $F = f^{\sigma}:Y \to Y$ where $\sigma = \min\{ n \ge1 : f^ny \in Y\}$.
Then $\mu=(\mu_0|Y)/\mu_0(Y)$ is an absolutely continuous invariant probability measure for $F$.
Define the induced roof function $\tau\to\R^+$ given by
$\tau(y) = \sum_{\ell=0}^{\sigma(y)-1} \tau_0(f^\ell y)$.
Let $F_t:Y^\tau\to Y^\tau$ be the corresponding suspension semiflow
with infinite invariant measure $\mu^\tau$.

Since $\tau_0$ is H\"older, it is standard that
$\mu(\tau>t)\sim c t^{-\beta}$ for some $c>0$ (see for example~\cite[Proposition~9.1]{BMTapp}).

\begin{prop} Suppose that $f_t$ has two periodic orbits (other than the neutral one) whose periods have irrational ratio.  Then hypotheses (H) and (S)(i) hold with
$\cB=\BV$ being the space of bounded variation functions on $Y$,
with norm $\|v\|_{\BV}=|v|_1+\Var v$.
\end{prop}

\begin{proof}
Hypotheses~(H)(i,iii) are verified in~\cite[Proposition~9.2]{BMTapp}.
Also, hypothesis~(H)(ii) is 
verified in~\cite[Proposition~9.2]{BMTapp} for $s\in\barH\cap B_\delta(0)$.

To complete the verification of~(H)(ii), we proceed as follows.
Since the density $d\mu/d\Leb$ lies in $\BV$ and is bounded above and below, it suffices to work with the non-normalised transfer operator
$\hat P(ib)v=P(e^{ib\tau}v)$ where 
$\int_Y Pv\,w\,d\Leb=\int_Y v\,w\circ F\,d\Leb$.

Let $\lambda=\inf g|_Y>1$.  Fix $L>0$.
It suffices to show that there exists a constant $C'$ such that
\[
\|\hat P(ib)^nv\|_{\BV}\le C'n|v|_1+C'n\lambda^{-n}\Var v,
\]
for all $|b|\le L$, $n\ge1$, $v\in \BV$.

Let $n\ge1$ and let $\{I\}$ be the partition of domains of branches for $F^n$.
There is a constant $C_0$ independent of $n$ such that
$\supI 1/(F^n)'\le C_0\diam I$ for all $I$.
Also $F'\ge\lambda$,
so $|1/(F^n)'|\le 1/\lambda^n$ for all $n$.

Write 
\[
\hat P(ib)^nv=\sum_I  \{\zeta_n\,e^{ib\tau_n}v\}\circ \psi_I\,1_{F^nI},
\]
where $\zeta_n=1/(F^n)'$, $\psi_I$ is the inverse branch $(F^n|_I)^{-1}$, and
$\tau_n=\sum_{j=0}^{n-1}\tau\circ F^j$ (not to be confused with $\tau_0$).
We have the standard estimate
\begin{align*}
|\hat P(ib)^nv|_1 & \le 
|\hat P(ib)^nv|_\infty\le \sum_I \supI (\zeta_n|v|)
\le \sum_I \supI \zeta_n  (\infI |v|+\Var_I v)
\\ & \le \sum_I  \supI \zeta_n (\diam I)^{-1}\int_I|v|
\;+\; \sum_I \lambda^{-n}\Var_I v
 \le C_0|v|_1+\lambda^{-n}\Var v.
\end{align*}

Next,
\begin{align*}
\Var & (\hat P(ib)^nv)  \le \sum_I  \Var_I(\zeta_n\,e^{ib\tau_n}v) + 2\sum_I \supI(\zeta_n|v|)
\\ & \le \sum_I \Var_I(\zeta_n v)
+\sum_I \supI(\zeta_n|v|)\Var_I e^{ib\tau_n}
+ 2C_0|v|_1+2\lambda^{-n}\Var v.
\end{align*}
A standard argument shows that 
\[
\sum_I \Var_I(\zeta_n v)\le C_1|v|_1+\lambda^{-n}\Var v,
\]
where $C_1=\sup_n |(F^n)''/[(F^n)']^2|$.
Also,
\[
\Var_I e^{ib\tau_n}\le |b|\Var_I\tau_n
\le L\sum_{j=0}^{n-1}\Var_I(\tau\circ F^j)
= L\sum_{j=0}^{n-1}\Var_{F^jI}\tau.
\]
Let $a$ be the domain of a branch for $F$.  Then 
$\tau|_a=\sum_{\ell=0}^{\sigma(a)-1}\tau_0\circ f^\ell$.
Since the images $f^\ell a$ are disjoint for $\ell<\sigma(a)$, it follows that
$\Var_a\tau\le \Var\tau_0$.  But $F^jI$ lies in such a domain $a$, so
$\Var_{F^jI}\tau\le \Var \tau_0$ and it follows that
$\Var_I e^{ib\tau_n}\le Ln\Var\tau_0$.
Hence
\[
\sum_I \supI(\zeta_n|v|)\Var_I e^{ib\tau_n}\le Ln\Var\tau_0
\sum_I \supI(\zeta_n|v|)
\le Ln\Var\tau_0\,(C_0|v|_1+\lambda^{-n}\Var v).
\]
Combining these estimates we have shown that
$\|\hat P(ib)^nv\|_{\BV}  \le (3C_0+C_1+C_0L\Var\tau_0)n|v|_1
+(4+L\Var\tau_0)n\lambda^{-n}\Var v$ as required.

Passing to the $L^2$ adjoint of $\hat R(ib)$, to verify (S)(i) it is equivalent to rule out the possibility that there exists 
$b\neq0$ and a  $\BV$ eigenfunction $v:Y\to S^1$ such that
$e^{ib\tau}v\circ F =v$.
Suppose that $y\in Y$ is a periodic point of period $k$ for $F$.
Now, $\BV$ functions have one-sided limits, and $F$ is orientation preserving, so
$v(y+)=v(F^k(y+))$.  
Substituting into the equation $e^{ib\tau_k}v\circ F^k=v$
we obtain $e^{ibq}=1$ where
$q=\tau_k(y+)$ is the period of the corresponding periodic orbit for $f_t$.
This is impossible under the periodic orbit assumption, so the $\BV$ eigenfunction $v$ cannot exist.
\end{proof}

It follows from Theorem~\ref{thm:H} that mixing for $F_t$
holds for all bounded almost everywhere continuous $\hat v$ supported in 
$H_n=\{(y,u)\in Y\times[0,\infty):\tau(y)-n\le u\le \tau(y)\}$ for some $n\ge1$,
and all $\hat w\in L^1(Y^\tau)$.

Let $v,w:[0,1]^{\tau_0}\to\R$ be observables where
$v$ is bounded and almost everywhere continuous and $w$ is integrable.
The projection $\pi:Y^\tau\to [0,1]^{\tau_0}$, $\pi(y,u)=f_u(y,0)$, defines a measure-preserving semiconjugacy from $F_t:Y^\tau\to Y^\tau$ to $f_t:[0,1]^{\tau_0}\to [0,1]^{\tau_0}$.
Define the lifted observables $\hat v=v\circ\pi$, $\hat w=w\circ\pi:Y^\tau\to\R$.
Then mixing for $F_t$ holds provided $\hat v$ is supported in an $H_n$ and hence the desired mixing result~\eqref{eq:ex} holds for $f_t$ and the observables $v$ and $w$.  
This includes all (finite linear combinations of) observables $v$ supported in $A\times [0,\inf_A\tau_0]$ where $A\subset \{\sigma\le j\}$ for some $j\ge1$.
(For such an observable $v$, we have $\supp\hat v\subset H_n$ for $n\ge j|\tau_0|_\infty$.)

We can enlarge the class of observables $v$ 
to include all bounded almost everywhere continuous functions that vanish on a neighborhood of the neutral fixed point.
First, by adjoining preimages of $Y$ we can enlarge $Y$ so that it contains $[\eps,1]$ for any prescribed $\eps>0$.
Hence we can suppose without loss that $\supp v\subset\{(x,u)\in [0,1]^{\tau_0}:x\in Y\}$.
Since $Y$ is the first return for $F$, it follows that
$\supp \hat v\subset\{(y,u)\in Y^\tau:u\le \tau_0(y)\}$.
Let $Y_j=\{y\in Y:\sigma(y)=j\}$.
Define 
$C'=C\times[0,|\tau_0|_\infty]$ where
$C=\bigcup_{j\ge1}\{Y_j:|\tau_0|_\infty<\inf_{Y_j}\tau\}$.
For the remaining $Y_j$, we have 
$j\le\inf_{Y_j}\tau\le|\tau_0|_\infty$ so
$\sup_{Y_j}\tau\le j|\tau_0|_\infty\le|\tau_0|_\infty^2$.
Hence 
$\supp\hat v\subset C'\cup H_n$ for $n\ge|\tau_0|_\infty^2$.
Such observables are covered by Section~\ref{sec:obs}: Take $\cC$ to be the collection of finite unions of intervals in $Y$ and
define $\cC'$ as in Proposition~\ref{prop:obs}.  Certainly $C'\in\cC'$.
Define $\cC''$ as in the proof of Theorem~\ref{thm:H}.  Then $C',H_n\in\cC''$,
so $C'\cup H_n$ lies in the ring generated by $\cC''$.  In particular, mixing holds for observables such as $\hat v$ supported in $C'\cup H_n$.

\begin{rmk}  To verify hypothesis (S)(ii) it suffices to rule out the possibility that there exists $b\neq0$, $\lambda\in S^1$ and a BV eigenfunction $v:Y\to S^1$ such that
$e^{ib\tau}v\circ F=\lambda v$.  But then every period $q=\tau(y)$ corresponding to a fixed point $y$ for $F$ satisfies
$e^{ibq}=\lambda$.
Hence hypothesis (S)(ii) holds provided this set of periods is not contained in a lattice of the form $a_1+a_2\Z$ for some $a_1,a_2>0$.
\end{rmk}

\begin{rmk} \label{rmk:PH2}
Combining this example with Remark~\ref{rmk:PH} leads to examples of partially hyperbolic intermittent flows preserving an infinite measure.  See~\cite{LiveraniTerhesiu16,M15} for similar examples in the discrete time invertible  setting.  In addition to extending to continuous time, our examples are an improvement over those in~\cite{LiveraniTerhesiu16,M15} as far as mixing is concerned, since we require no assumptions on smoothness of foliations (in contrast to~\cite{LiveraniTerhesiu16}) or Markov structure (in contrast to~\cite{M15}).
\end{rmk}

\subsection{Suspensions over unimodal maps}

Let $f_t:[0,1]^{\tau_0}\to[0,1]^{\tau_0}$ be a suspension over a unimodal map 
$f:[0,1]\to[0,1]$ as described in Example~\ref{ex:unimodal}.
We sketch the main ingredients following~\cite[Section~10]{BMTapp}.   

By~\cite[Lemma~10.2(a)]{BMTapp}, $\mu_0(\tau_0>t)=ct^{-\beta}+O(t^{-2\beta})$
where the constant $c>0$ is given explicitly.
By~\cite{Young98}, $f:[0,1]\to [0,1]$ is modelled by a Young tower $F:Y\to Y$ where $Y$ is a tower with exponential tails over a suitable inducing set $Z\subset [0,1]$.
The roof function $\tau_0$ lifts to a roof function $\tau:Y\to\R^+$ satisfying
$\mu(\tau>t)=ct^{-\beta}+O(t^{-2\beta})$ where
$\mu$ is the SRB measure on $Y$.  

To prove~\eqref{eq:ex}, it remains to verify hypotheses~(A)  and~(S)(i).
In~\cite[Section~8.1]{BMTapp}, a new function space $\cB$ is defined for Young towers with exponential tails, and 
hypothesis~(A)(i,ii) are verified.  This relies on a technical condition called (H3) in~\cite{BMTapp} which is verified in~\cite[Lemma~10.3]{BMTapp}.  (The Lasota-Yorke inequality (A)(ii) is proved in~\cite[Theorem~B.2]{BMTapp} for $s\in\barH\cap B_1(0)$ but holds equally for $s\in\barH\cap B_L(0)$ for any $L>0$.)
By~\cite[Proposition~8.6 and Lemma~10.4]{BMTapp}, hypothesis~(A)(iii) is satisfied.  Finally, hypothesis (S)(i) is immediate from the quasicompactness assumptions (A)(i,ii) and the assumption about periodic orbits for $f_t$.

\paragraph{Acknowledgements}
This research began as a result of discussions with Jon Aaronson, Henk Bruin, Dima Dolgopyat, P\'eter N\'andori, Fran\c{c}oise P\`ene and Doma Sz\'asz at
the thematic program \emph{Mixing Flows and Averaging Methods} at the
Erwin Schr\"odinger Institute (ESI), Vienna, April/May 2016.
We are particularly grateful to Bruin, Dolgopyat, N\'andori and P\`ene for continued discussions on this topic and to the referee for several helpful suggestions.

The research of DT and IM was supported in part by funding from ESI.
The research of IM was supported in part by a
European Advanced Grant {\em StochExtHomog} (ERC AdG 320977).

\def\polhk#1{\setbox0=\hbox{#1}{\ooalign{\hidewidth
  \lower1.5ex\hbox{`}\hidewidth\crcr\unhbox0}}}

\end{document}